\documentclass{amsart}
\usepackage[dvips]{graphicx}
\usepackage{amscd}
\usepackage{amsmath}
\usepackage{amsxtra}
\usepackage{amsfonts}
\usepackage{amssymb}

\newtheorem{theorem}{Theorem}[section]
\newtheorem{corollary}[theorem]{Corollary}
\newtheorem{lemma}[theorem]{Lemma}
\newtheorem{proposition}[theorem]{Proposition}
\newtheorem{conjecture}[theorem]{Conjecture}
\theoremstyle{definition}
\newtheorem{definition}[theorem]{Definition}
\newtheorem{remark}[theorem]{Remark}

\newtheorem{example}[theorem]{Example}
\theoremstyle{remark}

\renewcommand{\theclaim}{\textup{\theclaim}}

\newtheorem*{acknowledgements}{Acknowledgements}

\numberwithin{equation}{section}

\def\openone

{\mathchoice

{\hbox{\upshape \small1\kern-3.3pt\normalsize1}}

{\hbox{\upshape \small1\kern-3.3pt\normalsize1}}

{\hbox{\upshape \tiny1\kern-2.3pt\SMALL1}}

{\hbox{\upshape \Tiny1\kern-2pt\tiny1}}}

\makeatletter

\newbox\ipbox

\newcommand{\diracb}[1]{\left\langle #1\mathrel{\mathchoice

{\setbox\ipbox=\hbox{$\displaystyle \left\langle\mathstrut
#1\right.$}

\vrule height\ht\ipbox width0.25pt depth\dp\ipbox}

{\setbox\ipbox=\hbox{$\textstyle \left\langle\mathstrut
#1\right.$}

\vrule height\ht\ipbox width0.25pt depth\dp\ipbox}

{\setbox\ipbox=\hbox{$\scriptstyle \left\langle\mathstrut
#1\right.$}

\vrule height\ht\ipbox width0.25pt depth\dp\ipbox}

{\setbox\ipbox=\hbox{$\scriptscriptstyle \left\langle\mathstrut
#1\right.$}

\vrule height\ht\ipbox width0.25pt depth\dp\ipbox}

}\right. }

\newcommand{\dirack}[1]{\left. \mathrel{\mathchoice

{\setbox\ipbox=\hbox{$\displaystyle \left.\mathstrut
#1\right\rangle$}

\vrule height\ht\ipbox width0.25pt depth\dp\ipbox}

{\setbox\ipbox=\hbox{$\textstyle \left.\mathstrut
#1\right\rangle$}

\vrule height\ht\ipbox width0.25pt depth\dp\ipbox}

{\setbox\ipbox=\hbox{$\scriptstyle \left.\mathstrut
#1\right\rangle$}

\vrule height\ht\ipbox width0.25pt depth\dp\ipbox}

{\setbox\ipbox=\hbox{$\scriptscriptstyle \left.\mathstrut
#1\right\rangle$}

\vrule height\ht\ipbox width0.25pt depth\dp\ipbox}

} #1\right\rangle}

\newcommand{\bz}{\mathbb{Z}}
\newcommand{\br}{\mathbb{R}}

\newcommand{\bn}{\mathbb{N}}

\def\blfootnote{\xdef\@thefnmark{}\@footnotetext}

\renewcommand{\mod}{\operatorname{mod}}
\newcommand{\rtwomatrix}{\begin{bmatrix}2&0\\0&2\end{bmatrix}}
\newcommand{\rtwoonematrix}{\begin{bmatrix}2&1\\0&2\end{bmatrix}}
\newcommand{\rthreematrix}{\begin{bmatrix}3&0\\0&3\end{bmatrix}}
\newcommand{\rtwomatrixthreed}{\begin{bmatrix}2&0&0\\0&2&0\\0&0&2\end{bmatrix}}
\newcommand{\threevectoro}{\begin{bmatrix}0\\0\\0\end{bmatrix}}
\newcommand{\threevectorx}{\begin{bmatrix}1\\0\\0\end{bmatrix}}
\newcommand{\threevectory}{\begin{bmatrix}0\\1\\0\end{bmatrix}}
\newcommand{\threevectorz}{\begin{bmatrix}0\\0\\1\end{bmatrix}}
\newcommand{\threevectorxy}{\begin{bmatrix}1\\1\\0\end{bmatrix}}
\newcommand{\threevectorxz}{\begin{bmatrix}1\\0\\1\end{bmatrix}}
\newcommand{\threevectoryz}{\begin{bmatrix}0\\1\\1\end{bmatrix}}
\begin{document}
\title[Orthogonality in IFS]{Analysis of orthogonality and of orbits in affine iterated function systems}
\author{Dorin Ervin Dutkay}
\blfootnote{Research supported in part by the National Science
Foundation DMS 0457491}
\address[Dorin Ervin Dutkay]{\protect\raggedright
University of Central Florida,
	Department of Mathematics,
	4000 Central Florida Blvd.,
	P.O. Box 161364,
	Orlando, FL 32816-1364,USA}
\email{ddutkay@mail.ucf.edu}
\author{Palle E.T. Jorgensen}
\address[Palle E.T. Jorgensen]{Department of Mathematics, 
The University of Iowa, 14 MacLean Hall, Iowa City, IA 52242-1419, USA}
\email{jorgen@math.uiowa.edu}
\begin{abstract}
 We introduce a duality for Affine Iterated Function Systems (AIFS) which is naturally motivated by group duality in the context of traditional harmonic analysis. Our affine systems yield fractals defined by iteration of contractive affine mappings. We build a duality for such systems by scaling in two directions: fractals in the small by contractive iterations, and fractals in the large by recursion involving iteration of an expansive matrix. By a fractal in the small we mean a compact attractor $X$ supporting Hutchinson's canonical measure $\mu$ , and we ask when $\mu$ is a spectral measure, i.e., when the Hilbert space $L^2( \mu)$ has an orthonormal basis (ONB) of exponentials $\{e_\lambda\, |\, \lambda \in \Lambda\}$. We further introduce a Fourier duality using a matched pair of such affine systems. Using next certain extreme cycles, and positive powers of the expansive matrix we build fractals in the large which are modeled on lacunary Fourier series and which serve as spectra for $X$.  Our two main results offer simple geometric conditions allowing us to decide when the fractal in the large is a spectrum for $X$. Our results in turn are illustrated with concrete Sierpinski like fractals in dimensions 2 and 3.
\end{abstract}

\maketitle \tableofcontents
\section{Introduction}\label{intro}
While the world of fractals (see \cite{BD88}) entails both a fascinating geometry and analysis, the introduction of spectral theory into the subject tends to limit the number of possibilities, see e.g., \cite{AnLa97}. Intuitively and geometrically we think of a fractal as a set which ``looks the same'' at different scales, where scaling is defined relative to a family of transformations and the structures are studied under the name Iterated Function Systems (IFS), see e.g., \cite{Hut81}. 

       Varying the transformations then yields different classes of fractals. Motivated both by our problem and by our applications, we limit ourselves here to affine mappings. Iteration of these mappings then yields scales in the small and scales in the large.

       Two approaches to IFSs have been popular, one based on a discrete version of the more familiar and classical second order Laplace differential operator of potential theory, see \cite{KSW01, Kig04, LNRG96};  and the second approach is based on Fourier series, see e.g.,  \cite{JoPe98, DuJo05}. The first model in turn is motivated by infinite discrete network of resistors, and the harmonic functions are defined by minimizing a global measure of resistance, but this approach does not rely on Fourier series. In contrast, the second approach begins with Fourier series, and it has its classical origins in lacunary Fourier series \cite{Kah86}.

       Some of the more popular models for the potential theoretic approach center around concrete examples, and especially certain Sierpinski like fractals. These are various affine fractals which rely on a specific notion of self-similarity \cite{KiLa01}. In these fractals, the self-similarity is specified by a set of affine transformation in $d$-dimensional Euclidean space $\br^d$. This means that the fractals themselves, say $X$, are compact subsets of the ambient $\br^d$. While $X$ itself does not carry any linear structure, its ambient $\br^d$ does. Using a key idea of Hutchinson \cite{Hut81}, it is easy to see that every $X$ arises naturally as the support of an associated measure $\mu$, actually a family of measures. Consider the case when the family of mappings $(\tau_i)$ which define the IFS is finite, say $N$ maps, where each $\tau_i$ is contractive. As the maps are iterated, probabilities $(p_i)$ are assigned to the $N$ possibilities. Hence it is natural to ask when the Hilbert space $L^2( \mu)$ has an orthonormal basis (ONB) of exponentials $\{e_\lambda\, |\, \lambda \in \Lambda\}$ where $e_\lambda$ is $\exp( 2 \pi i\lambda \cdot x)$ restricted to $X$.  In that case $\mu$ is called a spectral measure, and the corresponding set $\Lambda$ is called a spectrum. Our first observation is that spectral measures must have equal probabilities, i.e., $p_i = 1/N$. As noted in (\ref{eqinv}) below, we restrict attention to this case. Motivated by examples, we further restrict to the case when the affine mappings $(\tau_i)$ are determined by a fixed invertible scaling matrix say $R$, and a finite set of translation vectors $B$ in $\br^d$.

        We introduce a duality for such Affine Iterated Function Systems (AIFS) which is naturally motivated by group duality in the context of traditional harmonic analysis, see e.g., \cite{HeRo70}. Nonetheless, our present objects $X$ are highly non-linear, and they are not groups. Since our affine systems are defined by iteration of invertible mappings, we rather think of them as fractals in the small and fractals in the large. By a fractal in the small we mean the compact attractor $X_B$ supporting the canonical measure $\mu_B$ of (\ref{eqinv}),  and we ask when $\mu_B$ is a spectral measure, i.e., when the Hilbert space $L^2( \mu_B)$ has an orthonormal basis (ONB) of exponentials $\{e_\lambda\, |\, \lambda \in \Lambda\}$?

     In the construction of this Fourier duality, a second system $(R^T, L)$ enters where $L$ is again a finite subset of $\br^d$ of the same cardinality as $B$. Using this set $L$, certain cycles called $W_B$-cycles, and positive powers of the transposed matrix $S = R^T$, we then proceed to build a fractal in the large $\Lambda = \Lambda(S,L)$.  Our main results Theorem \ref{thnoort} and Theorem \ref{propwb}, and Theorem \ref{thsier} offer simple geometric conditions allowing us to decide when $\Lambda(S,L)$ is a spectrum. Our results in turn are illustrated with concrete Sierpinski like fractals in dimensions 2 and 3. 

\par
            The Sierpinski examples fall in a subclass of AIFSs where the maps $\tau_i$ are similitudes. In our case, this is reflected in the scaling matrix $R$; it is a diagonal matrix. It is known that fractals $X$ built on similitudes have better separation properties, referring to the individual sets $\tau_i(X)$  (see \cite{BNR06}), and their Hausdorff dimension is known \cite{Fal97}.
\par
There are several versions of the planar Sierpinski examples. They were introduced originally (see \cite{Sier52} and  \cite{Ste95}) in the context of general topology, and in this context the generic topological properties are stable under most variations of the planar example. 

        In contrast, the role of the examples in spectral theory and in Fourier duality is of a later vintage (see e.g., \cite{LaWa97, JoPe98, StWa99, LaWa06}) and there the stability properties are quite different as can be seen from Section \ref{sier}.

        {\bf The Fuglede Problem.}  Our present analysis is motivated by what is known as ``the Fuglede Problem''; i.e., the problem of sorting out the relationship between Fourier spectrum and geometry for sets $\Omega$ in $\br^d$ of positive finite Lebesgue measure; see e.g., \cite{Fug74} and  \cite{Jor82}.  Actually, Fuglede \cite{Fug74} asked his question for bounded and open sets in $\br^d$. More precisely, here in our present paper we are asking the same question for fractals as the one Fuglede asked for ``classical'' domains in $\br^d$. Fuglede suggested that a measurable set $\Omega$ would allow an orthogonal basis ONB of complex exponentials $\{e_\lambda \,|\, \lambda \in \Lambda\}$ if and only if $\Omega$ is a tile for a single-tiling of $\br^d$ by translations, allowing overlap only on sets of Lebesgue measure zero. Each aspect of Fuglede's problem for domains in $\br^d$ is of independent interest; spectrum and geometry. While there already was up to 2004 a considerable amount of work on both sides of Fuglede's question, both on the spectral side (the {\it spectral sets}), and the tiling side, (see the references in \cite{JoPe98}) it was only recently \cite{Tao04} that the problem in its general form was shown to be negative. Tao's paper was quickly followed by several others.

                  Tao's counter example \cite{Tao04} to the Fuglede conjecture was only in one direction: Tao showed for $d = 5$ that there are spectral sets which do not tile ($\br^d$ by translations). The obstruction to tiling from Tao's example was a counting/divisibility argument which in fact motivated our present work on fractals. And there then came a counter example by Kolountzakis and Matolci ``in the other direction'' (a tile which is not a spectral set); specifically that sets which tile need not be spectral. 

This was then followed by Matolci's improvement of Tao's result down to four dimension, i.e., $d = 4$, 
\cite{Mat06}. An example of a non-spectral tile has now been claimed for $d=3$, \cite{FMM06}.

\section{Affine iterated function systems}\label{affi}
In this section we define the geometric structures under discussion, and we introduce our central themes: duality, measure, and orthogonality to be used later. In Definition \ref{defhada} we introduce the class of complex Hadamard matrices which link the two sides of our duality for affine iterated function systems (AIFSs.) 

       By their nature, iterated function systems (IFSs) give rise to combinatorial trees, to dynamics, and to associated cycles; see e.g., \cite{BD88, Fal97, Hut81}. They will be introduced as needed in our analysis below.

Let $R$ be a $d\times d$ integer matrix, which is expansive, i.e., all its eigenvalues $\lambda$ satisfy $|\lambda|>1$.
For a point $b\in\br^d$ we define the affine map
\begin{equation}\label{eqtaub}
\tau_b(x)=R^{-1}(x+b),\quad(x\in\br^d).
\end{equation}
Since $R$ is expansive, there exists a norm on $\br^d$ fo which $\|R^{-1}\|<1$.
\par
For a finite set $B\subset\br^d$ one can define the iterated function system $(\tau_b)_{b\in B}$.
\par
There exist a unique compact set $X_B$ with the property 
$$X_B=\bigcup_{b\in B}\tau_b(X_B).$$
The set $X_B$ is called the {\it attractor} of the IFS $(\tau_b)_{b\in B}$.
\par
Let $N$ be the cardinality of $B$.
There exists a unique probability measure $\mu_B$ on $\br^d$ such that:
\begin{equation}\label{eqinv}
\int f\,d\mu_B=\frac{1}{N}\sum_{b\in B}\int f\circ\tau_b\,d\mu_B,\quad(f\in C_c(\br^d)).
\end{equation}
The measure $\mu_B$ is called the {\it invariant measure} of the IFS $(\tau_b)_{b\in B}$.

\par
The Fourier transform of the measure $\mu_B$ is
$$\hat\mu_B(x)=\int e^{2\pi ix\cdot t}\,d\mu_B,\quad(x\in\br^d).$$
\par
Taking $f(t)=e^{2\pi ix\cdot t}$ in the invariance equation (\ref{eqinv}), one obtains
\begin{equation}\label{eqfinv}
\hat\mu_B(x)=m_B((R^T)^{-1}x)\hat\mu_B((R^T)^{-1}x),\quad(x\in\br^d),
\end{equation}

where
\begin{equation}\label{eqmb}
m_B(x)=\frac{1}{N}\sum_{b\in B}e^{2\pi ib\cdot x},\quad(x\in\br^d).
\end{equation}

Then, since $m_B(0)=1=\hat\mu_B(0)$ and since $m_B$ is Lipschitz and $(R^T)^{-1}$ is contractive, it follows that the following infinite product is absolutely convergent and
\begin{equation}\label{eqinfprod}
\hat\mu_B(x)=\prod_{n=1}^\infty m_B((R^T)^{-n}x),\quad(x\in\br^d).
\end{equation}

\begin{definition}\label{defhada}
Let $R$ and $B$ as above. For a finite subset $L$ of $\br^d$, we say that $(R,B,L)$ is a Hadamard triple if $L$ has the same cardinality as $B$, and the matrix
$$\frac{1}{\sqrt{N}}(e^{2\pi iR^{-1}b\cdot l})_{b\in B,l\in L}$$
is unitary.
\end{definition}

\begin{figure}
\setlength{\unitlength}{0.4\textwidth}
\begin{picture}(2.125,0.625)(0,-0.125)
\put(0,0){\includegraphics[bb=88 4 376 148,width=\unitlength]{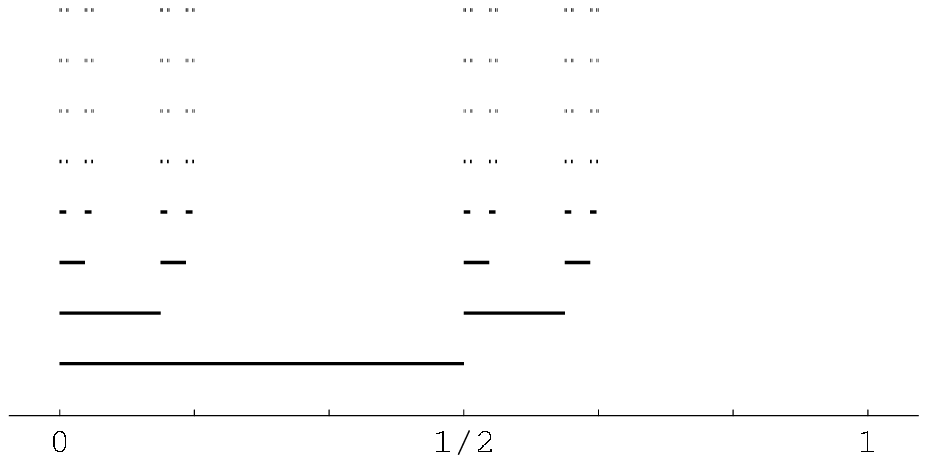}}
\put(1.125,0){\includegraphics[bb=88 4 376 148,width=\unitlength]{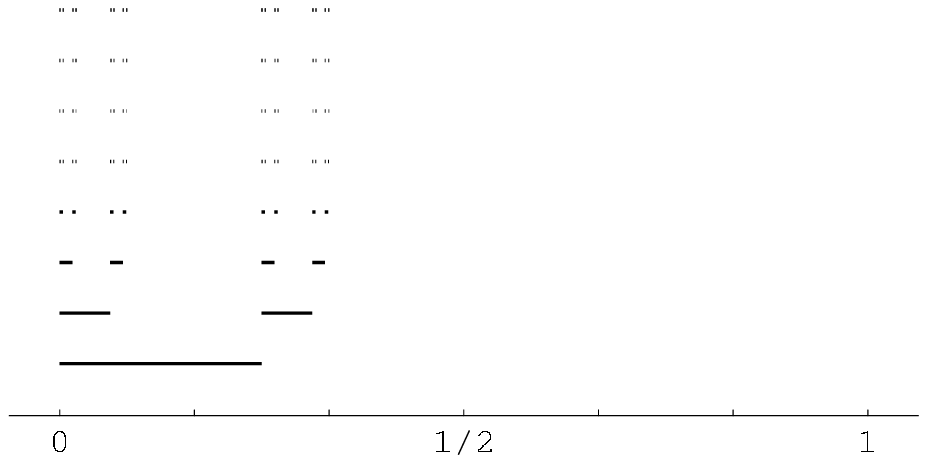}}
\put(0.5,-0.0625){\makebox(0,0)[t]{$X_B$}}
\put(1.625,-0.0625){\makebox(0,0)[t]{$X_L$}}
\end{picture}
\caption{The attractors $X_B$ and $X_L$ for $R = 4$, $B = \{0, 2\}$, $L = \{0, 1\}$}
\label{FigCantorFour}
\end{figure}

A simple example, with $d = 1$, of a system $(R, B, L)$ for which this unitarity
holds is: $R = 4$, $B = \{0, 2\}$, $L = \{0, 1\}$; see Figure \ref{FigCantorFour}, and \cite{JoPe98}. Both of
the attractors $X_B$ and $X_L$ has the same Hausdorff dimension 
$\log2/\log4$,
and
the measure $\mu_B$ is a spectral measure.
In Theorem \ref{thsier} below, we extend this result, not only to a complete analysis of the Sierpinski cases for $d = 1$, but also for $d = 2$ and $3$. So far, our results for $d = 4$ are only partial.

\begin{remark}\label{rem2diez}
Hadamard matrices with entries $\pm 1$ have a long history in combinatorics, see e.g., \cite{TaVu06}, while their complex variants with phase modulations $\exp(i 2\pi h)$ have a use in the study of communication filters. The complex variants are of a more recent vintage.

       Each complex Hadamard matrix has an associated matrix with real entries; its logarithmic variant, and its entries are the phase angles $h$. It is often referred to as the logarithmic form of the complex Hadamard matrix in question. Our condition in Definition \ref{defhada} above is that the numbers $h$ (the entries in the logarithmic variant) arise as inner products from two finite sets $B$ and $L$ of vectors in $\br^d$. These are the two sets which generate the respective sides, geometry and spectrum in our duality. 

       A useful fact about complex Hadamard matrices is that they are closed under tensor product, see e.g \cite{JoPe92}.

       The complex Hadamard matrices were introduced in spectral geometry, and in the study of spectral duality in the following papers \cite{Jor82, JoPe87, JoPe91, JoPe92, JoPe93}. Since then they have been used extensively in harmonic analysis, in varied contexts involving the analysis of Fourier bases and geometry in for example the papers \cite{Haa96, IoPe98, IoKaPe99,PeWa01, LaWa96, LaWa97, Tao04} among others. The original duality theme has expanded from its initial context to the harmonic analysis of fractals, the study of communications filters, of wavelets, and quantum theory.
\end{remark}

The function $W_B = | m_B|^2$ is the key to our duality consideration. Our aim is to study the spectral theory of the measure $\mu_B$. For this purpose the condition in Definition \ref{defhada} serves to identify a duality which holds for some but not for all affine iterated function systems (AIFS.)

       As in (\ref{eqtaub}) from the two pairs $(R, B)$ and $(S, L)$ we then define a dual pair of AIFSs, and we denote the respective compact attractors $X_B$ and $X_L$. Here we set $S = R^T$ (the transposed matrix.)  

       Our use of duality for such pairs is motivated by group duality. Nonetheless, we stress that our objects are not groups. Rather we think of them as fractals in the small and fractals in the large. By a fractal in the small we mean the compact attractor $X_B$ supporting the canonical measures $\mu_B$ of (\ref{eqinv}),  and we ask when $\mu_B$ is a spectral measure, i.e., when the Hilbert space $L^2( \mu_B)$ has an orthonormal basis (ONB) of exponentials  $\{e_\lambda\, |\, \lambda \in \Lambda\}$? Here $e_\lambda(x)=e^{2\pi i\lambda\cdot x}$. If $\mu_B$ is a spectral measure, there is such a set $\Lambda$, and we call $\Lambda$ a spectrum.

     In understanding this construction, the second system $(S, L)$ enters. Using $L$, certain cycles called $W_B$-cycles, and positive powers of $S$ we then proceed to build a fractal in the large $\Lambda(S,L)$, and the aim is to decide when $\Lambda(S,L)$ is a spectrum or an ONB for $L^2(\mu_B)$. Associated to $B$, we define $m_B$ as in (\ref{eqmb}), and the weight function $W_B = | m_B |^2$. The condition in Definition \ref{defhada} implies that $W_B$ satisfies the following normalization

\begin{equation}\label{eqnorm}
\sum_{l\in L}W_B(\tau_lx)=1,\quad(x\in\br^d),
\end{equation}
where 
$$\tau_l(x)=S^{-1}(x+l),\quad(x\in\br^d,l\in L).$$

       Our AIFS-fractals and their spectra are mathematical counterparts of a theme in solid state physics; see [Sen95]: Atoms in quasicrystals are arranged in a definite and orderly way, but it is not periodic. The periods are ``disturbed'' by $N$-point diffraction; in physics called X-ray diffraction. This spectral theoretic view of diffraction too involves generalized Fourier expansions going beyond the original and more familiar periodic case based on lattices in $\br^d$, \cite{HeRo70}. The book \cite{Sen95} is a delightful exposition covering such a variety of quasiperiodic geometries, and starting with those first observed in diffraction theory of quasicrystals from solid state physics. Senechal's book \cite{Sen95} further contains complete details and references to the research literature.

     While our aim here is quite different from that of diffraction theory, they both involve a certain Fourier spectrum based on distributions of sample points in $\br^d$. For example our function $W_B = |m_B|^2$ occurs for related finite sets $B$, and a variant of it explains $N$-point diffraction in physics, see \cite[Ch 3, especially p. 83]{Sen95}. 

     The probabilistic significance of our function $W_B$ is contained in (\ref{eqnorm}). Here we use $W_B$ in order to define transition probabilities: If $x$ and $y$ are points in $X_L$ such that $\tau_l(y) = x$, for some $l \in L$, then $W_B(y)$ is the probability of transition from $x$ to $y$, see \cite{DuJo05}. In fact, this random-walk approach relates the geometry of the initial fractal $X_B$ (in the small) to spectral data in the large, where the computation of spectrum again uses iteration of a finite dual system of affine maps in the ambient space $\br^d$.

\section{Non-spectral measures}\label{non}
   In this section we answer the following question:  When a measure $\mu$ arises from an affine fractal iteration taking place in $\br^d$, then what is the largest number of orthogonal complex exponentials in $L^2(\mu)$ ?; i.e., what is the cardinality of the largest orthogonal $\{e_\lambda\}$ family in $L^2(\mu)$?

        An earlier theorem \cite{JoPe98} for $d = 1$ states that if $X_3$ is the classical triadic Cantor set of fractal dimension $s = \log_32$ with associated fractal measure $\mu_3$ on the line, then $\mu_3$  is not spectral. In fact, there are no more than two orthogonal $e_\lambda$'s in $L^2(\mu_3)$. This is in contrast to $\mu_3$'s natural counterpart $\mu_4$, Cantor's construction in base four of fractal dimension $s =1/2$ where (by \cite{JoPe98}) $L^2(\mu_4)$ even has a whole ONB of exponentials $e_\lambda$. This section is concerned with the shades in between for affine fractals in $\br^d$.

Starting with some given AIFS  $(R, B, \mu_B )$  in $\br^d$, with the measure $\mu_B$ depending on both the fixed matrix $R$, and the subset $B$, our two main questions are to estimate the number of orthogonal complex exponentials in $L^2(\mu_B)$, and to find them. There might in fact be no more than two orthogonal complex exponentials. Generally for AIFSs, the possibilities fall in one of the following three classes: 
\begin{itemize}
\item[(a)] There can be at most a finite number of orthogonal complex exponentials in $L^2(\mu_B)$.
\item[(b)] There are natural infinite families of orthogonal complex exponentials.  
\item[(c)] One of the infinite families from (b) is in fact an orthonormal basis (ONB) in $L^2(\mu_B)$. If (c) holds, we say that $\mu_B$ is a spectral measure.
\end{itemize}

      Our first result Theorem \ref{thnoort} gives conditions for when some  AIFS  $(R, B, \mu_B )$  falls in class (a). It is phrased in terms of certain orbits defined from the transposed matrix $S = R^T$ and the set $B$.

\begin{theorem}\label{thnoort}
Let $R$ be a $d\times d$ expansive integer matrix, and $B\subset\bz^d$ some finite set of cardinality $N$. Consider the IFS $\tau_b(x)=R^{-1}(x+b)$, $b\in B, x\in\br^d$, and let $\mu_B$ be its invariant measure. Let
$$m_B(x)=\frac{1}{N}\sum_{b\in B}e^{2\pi ib\cdot x},\quad(x\in\br^d).$$
Let $Z$ be the set of the zeroes of $m_B$ in $[0,1)^d$.
\begin{enumerate}
\item
If $Z$ is contained in a set $Z'\subset[0,1)^d$ of finite cardinality $|Z'|$, which does not contain $0$, and that satisfies the property 
$$R^T(Z'+\bz^d)\subset Z'+\bz^d.$$
Then there exist at most $|Z'|+1$ mutually orthogonal exponential functions in $L^2(\mu_B)$. In particular, $\mu_B$ is not a spectral measure. 
\par
The set $Z'$ can be taken to be $$Z':=\{(R^T)^k x\mod\bz^d\,|\,k\geq0, x\in[0,1)^d, m_B(x)=0\}$$ if this set is finite and does not contain $0$. 

\item
Let $\mathcal{O}(Z)$ be the orbit of $Z$ under the map $x\mapsto R^Tx\mod\bz^d$, i.e.,
$$\mathcal{O}(Z):=\{(R^T)^nx\mod\bz^d\,|\,x\in Z,n\in\bn\}.$$
Assume that the Euclidian distance from $\mathcal{O}(Z)$ to $\bz^d$ is $\delta>0$. Then there exists at most $\left(\left\lfloor \frac{\sqrt{d}}{\delta}\right\rfloor+1\right)^d$ mutually orthogonal exponential functions in $L^2(\mu_B)$.
\end{enumerate}
\end{theorem}

\begin{proof} (i)
Suppose there exists a family of mutually orthogonal exponential functions $\{e_{\lambda}\,|\,\lambda\in\Lambda\}$, with $|\Lambda|>|Z'|+1$. By taking some $\lambda_0\in\Lambda$ and replacing $\Lambda$ by $\Lambda-\lambda_0$, we may assume that $0\in\Lambda$.
\par
The orthogonality implies that, for $\lambda,\lambda'\in\Lambda$, with $\lambda\neq\lambda'$, one has 
\begin{equation}\label{eqzmub}
\hat\mu_B(\lambda-\lambda')=0.
\end{equation}
{}From the infinite product formula (\ref{eqinfprod}) $\mu_B$, we obtain that for some $n\geq1$, $(R^T)^{-n}(\lambda-\lambda')$ is a zero for $m_B$ so it is in $Z'+\bz^d$, because $B\subset\bz^d$ so $m_B$ is $\bz^d$-periodic.
Using the hypothesis we get that 
\begin{equation}\label{eqlala}
\lambda-\lambda'\in Z'+\bz^d.
\end{equation}
\par
Let $M:=\{\lambda\mod [0,1)^d\,|\,\lambda\in\Lambda\}$. We claim that $|M|\leq |Z'|+1$.
\par
Indeed, we have 
$$M+\bz^d\subset \Lambda+\bz^d\subset (\{0\}\cup Z')+\bz^d,$$
where we used (\ref{eqlala}) with $\lambda'=0$ for the last inclusion.
Thus $M\subset\{0\}\cup Z'$, and the cardinality of $M$ is $p\leq |Z'|+1$.
\par
We can enumerate the elements $M=\{m_1,\dots,m_p\}$ and 
$$\Lambda\subset\cup_{k=1}^p(m_k+\bz^d).$$
But since $|\Lambda|\geq |Z'|+1\geq p$, one of the sets $m_k+\bz^d$ will contain two distinct elements in $\Lambda$. Hence, 
there exist $\lambda\neq\lambda'\in\Lambda$ such that $\lambda-\lambda'\in\bz^d$. But this will contradict (\ref{eqlala}), because $0\not\in Z'$.
\par
(ii) As in (i) we can assume that the family of mutually orthogonal functions $\Lambda$, contains $0$. The same argument as in (i) shows that 
\begin{equation}\label{eqoz}
\lambda-\lambda'\in\mathcal{O}(Z)+\bz^d,\quad(\lambda\neq\lambda', \lambda,\lambda'\in\Lambda).
\end{equation}
Suppose $\Lambda$ contains more than $k^d$ functions where $k:=\left\lfloor \frac{\sqrt{d}}{\delta}\right\rfloor+1$. Then we can divide the cube $[0,1)^d$ into $k^d$ cubes with sides equal to $\frac1k$. One of these small cubes will contain two elements $\lambda\mod\bz^d$ , $\lambda'\mod\bz^d$, with $\lambda,\lambda'\in\Lambda$ and $\lambda\neq\lambda'$. But then the Euclidian distance between $\lambda\mod\bz^d$ and $\lambda'\mod\bz^d$ is less than the diagonal which is 
$$\frac{\sqrt{d}}k=\frac{\sqrt{d}}{\left\lfloor \frac{\sqrt{d}}{\delta}\right\rfloor+1}<\delta.$$
Therefore 
$$d((\lambda-\lambda')\mod\bz^d,\bz^d)=d(\lambda\mod\bz^d-\lambda'\mod\bz^d,\bz^d)\leq$$$$ d(\lambda\mod\bz^d-\lambda'\mod\bz^d,0)<\delta.$$
But, from (\ref{eqoz}) we have that $(\lambda-\lambda')\mod\bz^d$ is in $\mathcal{O}(Z)$ hence, its distance to $\bz^d$ is at least $\delta$. This contradiction yields (ii).
\end{proof}

\begin{example}
Take 
$$R:=\left[\begin{array}{cc}2&1\\0&2\end{array}\right],\quad B:=\{(0,0),(1,0),(0,1)\}.$$
Then $\mu_B$ is not a spectral measure.
\par
We have $$m_B(x,y)=\frac{1}{3}(1+e^{2\pi ix}+e^{2\pi iy}),\quad(x,y\in\br).$$
As beforein the proof of Theorem \ref{thsier}(ii), we get that the zeroes of $m_B$ inside $[0,1)^2$ are $(1/3,2/3)$ and $(2/3,1/3)$. Iterating the map $x\mapsto (R^T)x\mod\bz^2$ on these points, we see obtain the following cycle
$$(1/3,2/3)\rightarrow(2/3,2/3)\rightarrow(1/3,0)\rightarrow(2/3,1/3)\rightarrow(1/3,1/3)\rightarrow(2/3,0)\rightarrow(1/3,2/3).$$
Take $Z'$ to be the set of these $6$ points. Applying Theorem \ref{thnoort} we see that $\mu_B$ cannot have more than $7$ mutually orthogonal exponential functions.
\end{example}
\begin{example}\label{exp2}
Take $R=2$ and $B=\{0,a\}$, with $a\in\br$, $a\neq 0$, but assign different probabilities $p_1\neq p_2$ to $\tau_0$, $\tau_a$, respectively; $p_1+p_2=1$. The invariance equation (\ref{eqinv}) becomes
$$\int f\,d\mu_{B,p}=p_1\int f\circ\tau_0\,d\mu_{B,p}+p_2\int f\circ\tau_a\,d\mu_{B,p},$$
and (\ref{eqinfprod}) is true with 
$$m_B(x):=m_{B,p}(x):=p_1+p_2e^{2\pi iax},\quad(x\in\br).$$
We claim that no two exponential functions are orthogonal in $L^2(\mu_{B,p})$. 
\par
For this, note first that $m_{B,p}(x)=0$ implies $-p_1=p_2e^{2\pi iax}$ so $|p_1/p_2|=1$, which we assumed not to be true. So $m_{B,p}$ has no zeroes, and therefore, with (\ref{eqinfprod}), $\hat\mu_{B,p}$ has no zeroes, and this proves our claim.
\end{example}

\section{Duality}\label{dual}
   Here we extend the duality notions which originated in The Fuglede Problem, see section \ref{intro}, the relationship between ONBs of complex exponentials $e_\lambda$ and geometry. For our fractal measures $\mu$ (with support on fractal Sierpinski sets in $\br^d$) we now ask when $\mu$ allows an orthogonal basis ONB of complex exponentials $\{e_\lambda\, |\, \lambda \in \Lambda\}$. But of course since we are now dealing with iteration limits, also called affine IFSs,  the modified ONB property refers to the Hilbert space $L^2(\mu)$, and not to the restriction of the $d$-dimensional Lebesgue measure. Yet, we are still considering restrictions to the affine IFS-fractal of the $e_\lambda$ functions on $\br^d$. Note that depending on $d$, our Sierpinski-Hutchinson measures $\mu$ have fractal dimension of value $s$ where $s$ is typically smaller than $d$. We then modify the other side of the Fuglede problem, the geometric side, and address the relation between spectrum and geometry. But it is different because the fractals are different creatures, for example they do not have a linear structure. The restricted class we consider is instead determined by an expansive $d$ by $d$ matrix $R$ and a finite and fixed set $B$ of given vectors in $\br^d$. By iterating powers of $R$, we arrive at fractals in the small and what may be called fractals in the large. It is the latter that typically may serve as spectrum of the former. Because of Fourier duality we must work with the transposed scaling matrix, i.e., $R^t$ on the ONB-dual side. Our main theorems concern the detailed dependence on $d$ of the answer to the spectral questions.

 Starting with a given AIFS  $(R, B, \mu_B )$  in $\br^d$, we introduced a second system $(S, L )$ dual to the first, where $S = R^T$, and where the dual subset $L$ in $\br^d$ is chosen according to Definition \ref{defhada}. The attractor for the first system is denoted $X_B$, and for the second $X_L$. We will resume our study of orthogonal complex exponentials in $L^2(\mu_B)$.  In Theorem \ref{propwb} we show that the possibilities for orthogonal complex exponentials depend on a certain class of cycles in $X_L$ called $W_B$-cycles, and we will outline the interconnections.

A cycle $C$ for the $(S,  L)$-system which has the additional property that $W_B(x) = 1$ for all $x$ in $C$ is called a $W_B$-cycle.

     As illustrated in Section \ref{sier} for the class of Sierpinski fractals, the $W_B$-cycles play an important role in our understanding of those affine measures $\mu_B$ which are also spectral measures, i.e., measures $\mu$ for which $L^2( \mu)$ has an orthonormal basis (ONB) of exponentials $\{e_\lambda\, | \lambda \in \Lambda\}$.

In this section we show more generally how the $W_B$-cycles may be accounted for by a certain lattice structure.

\begin{theorem}\label{propwb}
Suppose the exist $d$ linearly independent vectors in the lattice $\Gamma$ generated by $$\{\sum_{k=0}^nR^kb_k\,|\,b_k\in B,n\in\bn\}.$$ Define
$$\Gamma^\circ:=\left\{x\in\br^d\,|\,\beta\cdot x\in\bz\mbox{ for all }\beta\in \Gamma\right\}.$$
Then $\Gamma^\circ$ is a lattice that contains $\bz^d$, is invariant under $S$, and if $l,l'\in L$ with $l-l'\in S\Gamma^\circ$ then $l=l'$. Moreover 
$$\Gamma^\circ\cap X_L=\cup\{C\,|\, C\mbox{ is a }W_B\mbox{-cycle}\}.$$ 
\end{theorem}
\begin{proof}
The fact that $\Gamma^\circ$ is a discrete lattice follows from the existence of the $d$ linearly independent vectors. Since $0\in B$, it follows that $S\Gamma^\circ\subset\Gamma^\circ$. Since $B\subset\bz^d$ and $R$ has integer entries, we get that $\Gamma^\circ$ contains $\bz^d$.
\par
Take now some $W_B$-cyle, $C:=\{x_0,x_1:=\tau_{l_1}x_0,\dots,x_{p-1}:=\tau_{l_{p-1}}\cdots\tau_{l_1}x_0\}$, with $\tau_{l_p}\cdots\tau_{l_1}x_0=x_0$. We have $Sx_0\equiv x_{p-1}\mod\bz^d$, and by induction $S^kx_0$ is congruent to some point in $C$ modulo $\bz^d$. Therefore $W_B(S^kx_0)=1$ for all $k\geq0$.
\par
But then $W_B(x_0)W_B(Sx_0)\cdots W_B(S^{n}x_0)=1$ so 
$$\left|\sum_{b_0,\dots,b_n\in B}\prod_{k=0}^ne^{2\pi i b_k\cdot S^kx_0}\right|=N^{n+1}.$$
Since there are $N^{n+1}$ terms in the sum, and all of them have modulus $1$, it follows that all the terms are equal to $1$. Therefore, for all $b_0,\dots,b_n\in B$,
$$e^{2\pi i(\sum_{k=0}^nR^kb_k)\cdot x}=1$$
which implies that $x\in\Gamma^\circ$.
\par
Conversely, take $x_0\in X_L\cap\Gamma^\circ$. Then $x_0=\sum_{k=1}^\infty S^{-k}l_k$ for some $l_k\in L$. 
Define the points $x_p$ obtained from shifting the expansion of $x_0$:
$$x_p:=\sum_{k=1}^\infty S^{-k}l_{k+p},\quad(p\geq0).$$
Then note that $x_{p+1}=Sx_p-l_{p+1}$. By induction, since $x_0\in\Gamma^\circ$ and $\Gamma^\circ$ is invariant under $S$ and contains $\bz^d$, it follows that $x_{p}\in\Gamma^\circ$, for all $p\geq0$. Also $x_p\in X_L$ so $x_{p}\in\Gamma^\circ\cap X_L$ for all $p\geq 1$. However , since $\Gamma^\circ$ is discrete and $X_L$ is compact, this intersection is finite, so there exist 
$p\geq 0$ and $m\geq 1$ such that $x_p=x_{p+m}$.
\par
We claim that if $l,l'\in L$ and $l\neq l'$ then $l-l'\not\in S\Gamma^\circ$. Suppose not. Then
$S^{-1}(l-l')\cdot b\in\bz$ for all $b\in B$. Then
$$\sum_{b\in B}e^{2\pi ib\cdot S^{-1}(l-l')}=N,$$
and this contradicts the Hadamard property.
\par
We have $x_p=S^px_0-S^{p-1}l_1-S^{p-2}l_2-\dots-l_p$ and $x_p=x_{p+m}=S^{p+m}x_0-S^{p+m-1}l_1-\dots-l_{p+m}$.
\par
Using the previous claim, since $x_0\in\Gamma^\circ$, and $L\subset\Gamma^\circ$, we get $l_{p+m}=l_{p}$, so $x_{p-1}=S^{-1}(x_{p}+l_p)=S^{-1}(x_{p+m}+l_{p+m})=x_{p+m-1}$.
By induction $x_0=x_{m}$, and since $x_0=\tau_{l_1}\cdots\tau_{l_{m-1}}x_m$, this proves that $\{x_0,x_1,\dots,x_{m-1}\}$ is a cycle inside $\Gamma^\circ\cap X_L$.
\par
Since for all $x\in\Gamma^\circ$, one has $b\cdot x\in\bz$ for all $b\in B$, it follows that $W_B(x)=1$ for points in $\Gamma^\circ$. Hence $\{x_0,x_1,\dots,x_{m-1}\}$ is a $W_B$-cycle.

\end{proof}

The following lemma shows that our conclusions about the spectral properties of $\mu_B$ are invariant under a linear change of coordinates in $\br^d$. It applies for example to Variation 2 above, and it is used again at the end of Section \ref{dual}.
\begin{lemma}\label{lemtran}
Let $V$ be a $d\times d$ (real) matrix. Let $R_V:=VRV^{-1}$, $B_V:=VB$. Consider the invariant measures $\mu_B$ and $\mu_V:=\mu_{B_V,R_V}$. Then
$$\hat\mu_V(x)=\hat\mu_B(V^Tx),\quad(x\in\br^d).$$
$\mu_B$ is spectral iff $\mu_V$ is spectral. Suppose $\mu$ is spectral with spectrum $\Lambda$, and set $\Lambda_V:=(V^T)^{-1}\Lambda$. Then $\Lambda_V$ is a spectrum for $\mu_V$.
\end{lemma}
\begin{proof}
We have
$$\hat\mu_V(x)=\prod_{n\geq1}m_{VB}((R_V^T)^{-n}x)=\prod_{n\geq1}m_B((R^T)^{-n}V^Tx)=\hat\mu_B(V^Tx).$$
$\mu_B$ is spectral iff $\mu_V$ is spectral. Suppose $\mu$ is spectral with spectrum $\Lambda$, and set $\Lambda_V:=(V^T)^{-1}\Lambda$. Then $\Lambda_V$ is a spectrum for $\mu_V$.
\par
{}From this equality it follows that for continuous compactly supported functions on $\br^d$, 
$$\int f(x)\,d\mu_V(x)=\int f(Vx)\,d\mu_B(x).$$
Therefore the map $\Phi(f)(x)=f(Vx)$ is an isometry between $L^2(\mu_V)$ and $L^2(\mu_B)$, and 
$\Psi(e_\lambda)=e_{V^T\lambda}$. This implies the last statement. 
\end{proof}

\section{Sierpinski fractals}\label{sier}
       In this section we apply our main results to a special class of AIFS  $(R, B, \mu_B )$  in $\br^d$, called Sierpinski fractals. What sets them apart from the other AIFSs is our choice of the set $B$ of translation vectors, $B = \{0, e_1,..., e_d\}$ where $e_i$ denote the $i$-th canonical basis vector in $\br^d$. Even though the best known Sierpinski fractal \cite{Sier52} was initially only envisioned for $d = 2$, it seems natural to refer to the entire class as the Sierpinski fractals.

       We continue our focus on the class of Sierpinski examples, but we increase the dimension of the ambient space, i.e., the $\br^d$  containing the vertices $B$ of our Sierpinski attractor $X_B$, and we outline the changes in the conclusions above from 2D regarding spectrum, scaling rules, and orthogonality relations. There are several reasons why the class of Sierpinski examples is of independent interest. As already noted, it is widely studied; but in addition, a recent paper of J. d'Andrea, K. Merrill, and J. Packer \cite{DMP06} shows that the Sierpinski structures play a key role in a certain image processing algorithm. This fact is also based in part on a certain design (due to the co-authors) of a multiresolution analysis for the present affine fractal systems (AIFS); see \cite{DuJo06b}. With our introduction here of scaling in the small and in the large, our present paper has in common with \cite{DMP06} the use of nested scales spaces and of recursive algorithms.

\begin{figure}
\setlength{\unitlength}{0.32\textwidth}
\begin{picture}(3.125,3.875)(0,-0.125)
\put(0.05,2.75){\includegraphics[bb=88 4 376 292,width=\unitlength]{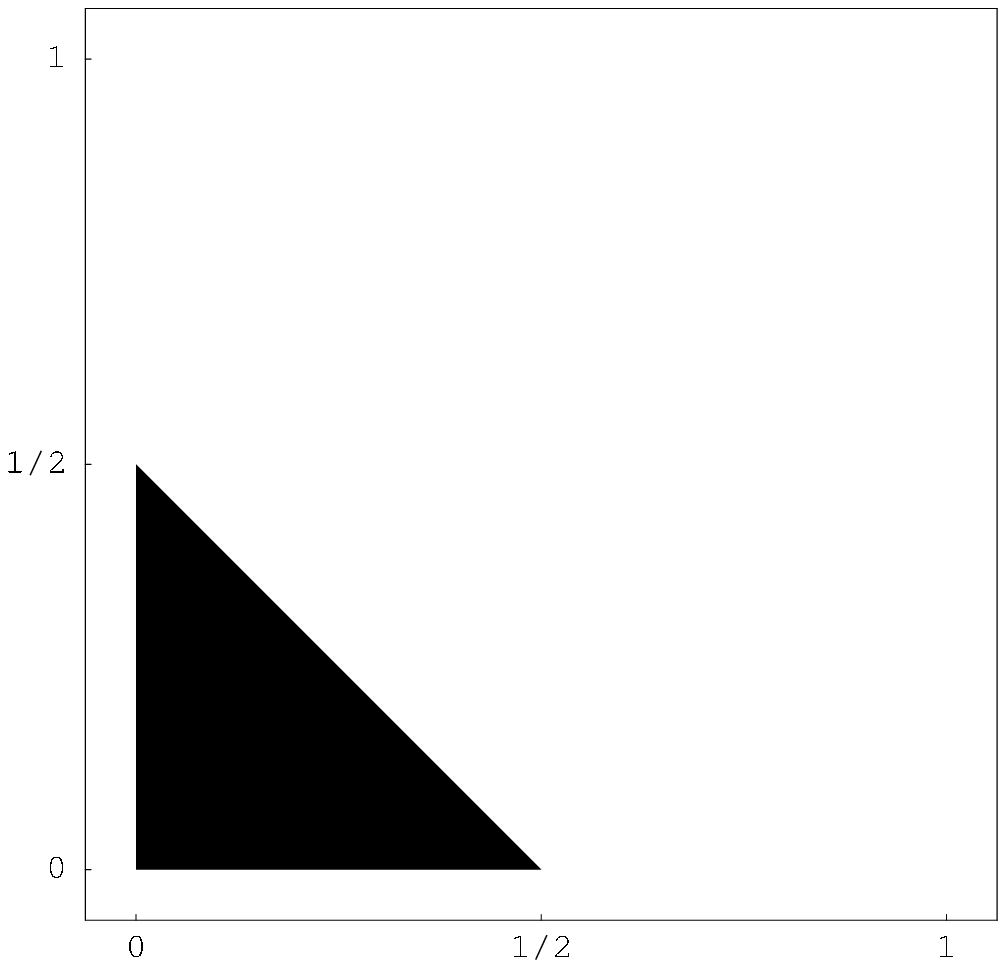}}
\put(1.0875,2.75){\includegraphics[bb=88 4 376 292,width=\unitlength]{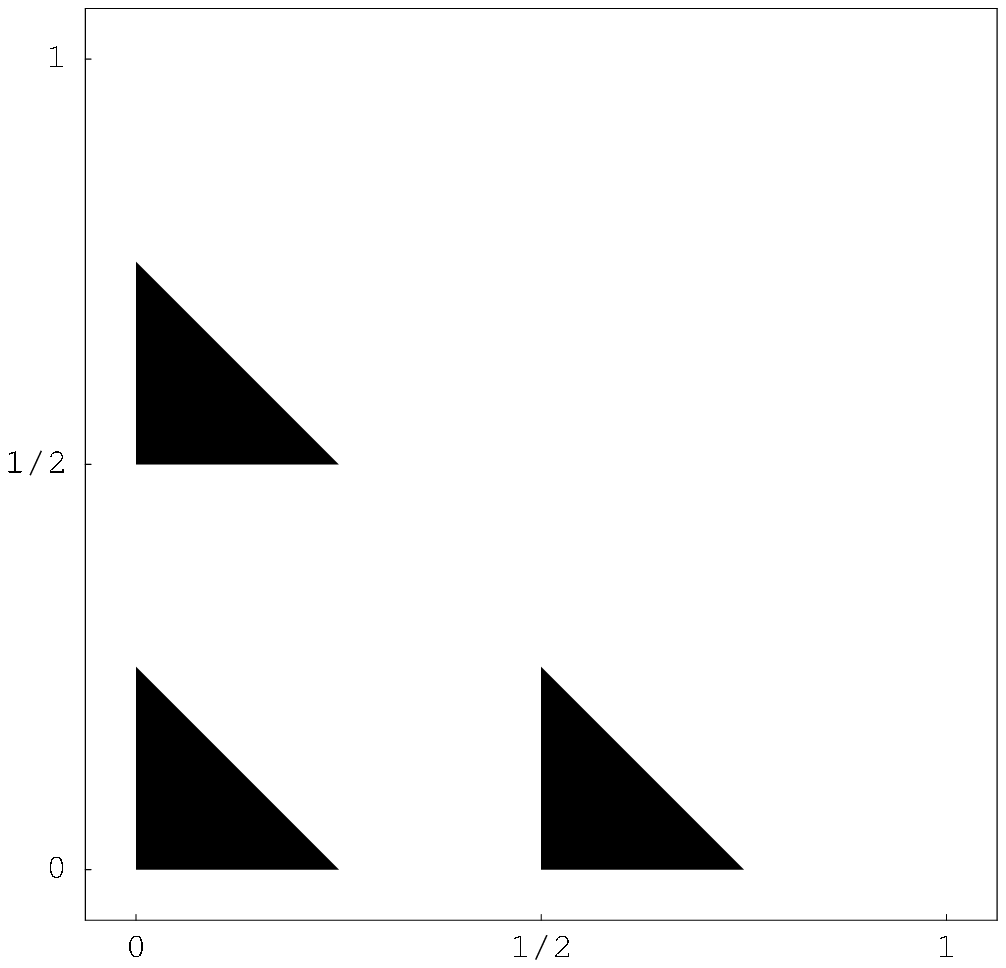}}
\put(2.125,2.75){\includegraphics[bb=88 4 376 292,width=\unitlength]{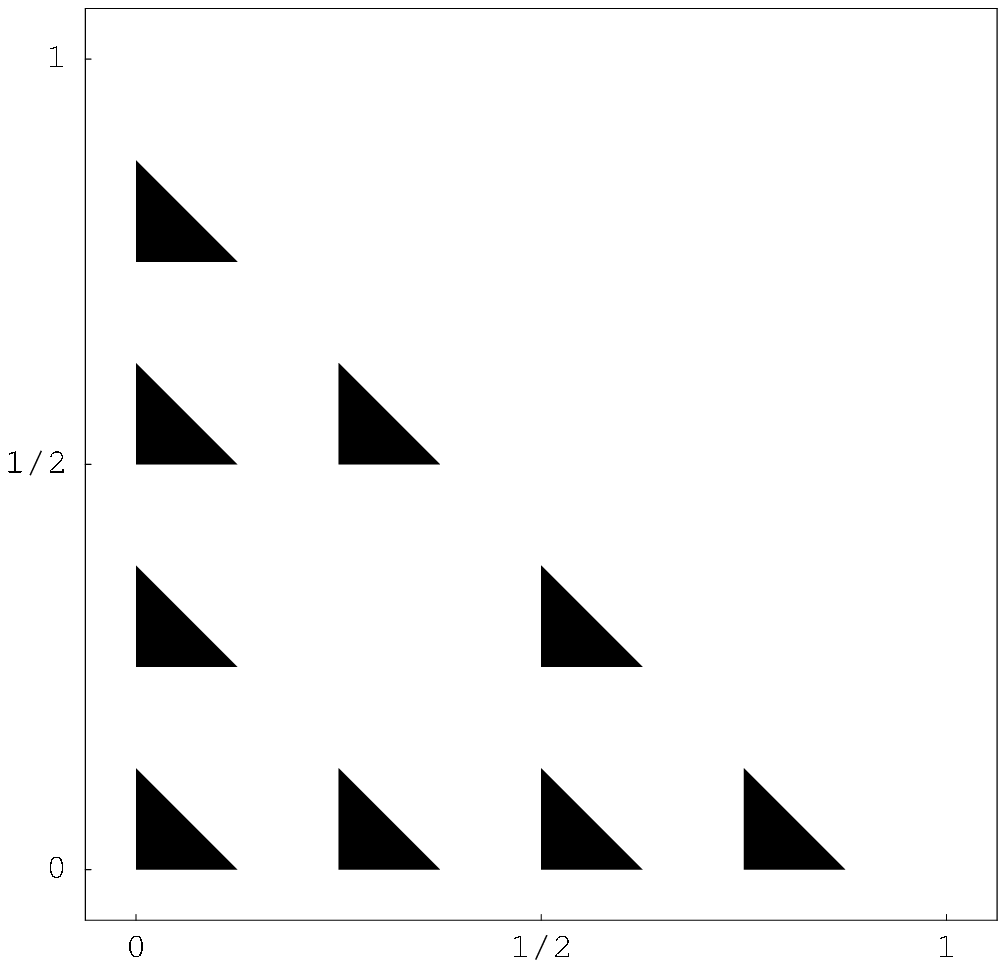}}
\put(0.60,2.75){\makebox(0,0)[t]{First iteration}}
\put(1.6375,2.75){\makebox(0,0)[t]{Second iteration}}
\put(2.675,2.75){\makebox(0,0)[t]{Third iteration}}
\put(0.05,1.625){\includegraphics[bb=88 4 376 292,width=\unitlength]{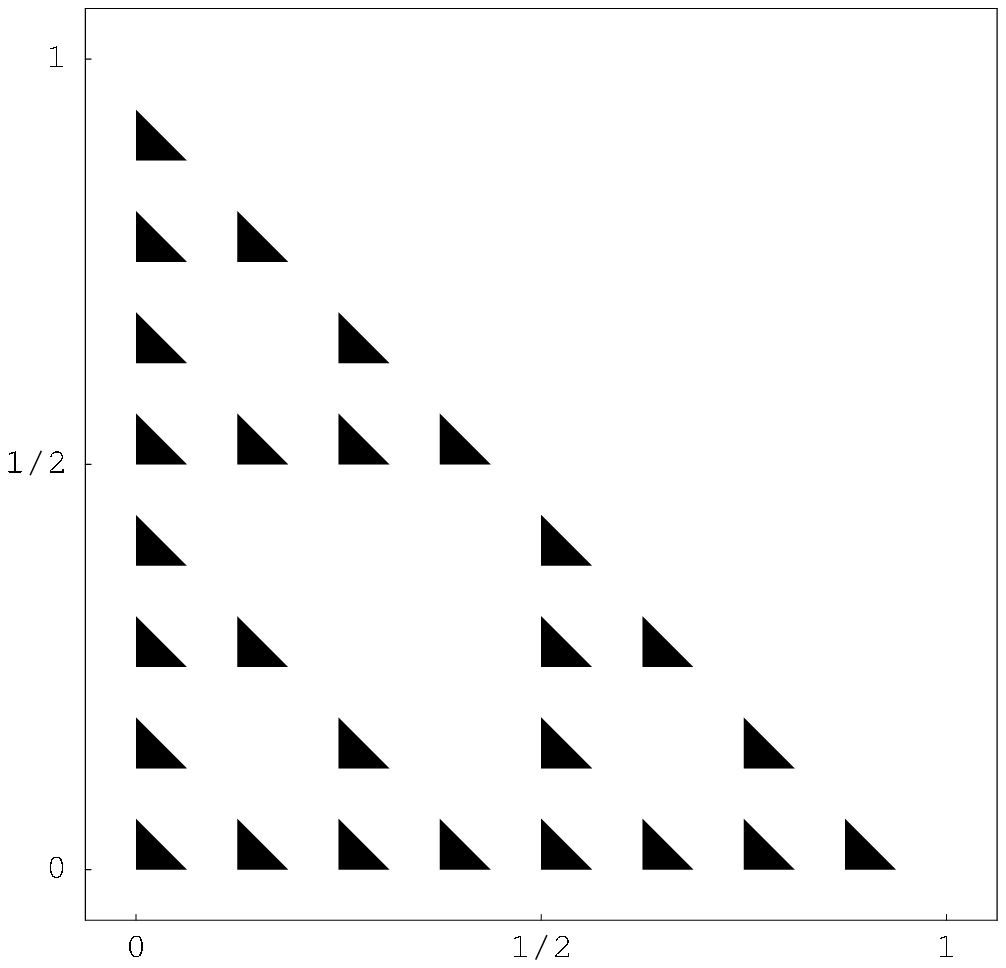}}
\put(1.0875,1.625){\includegraphics[bb=88 4 376 292,width=\unitlength]{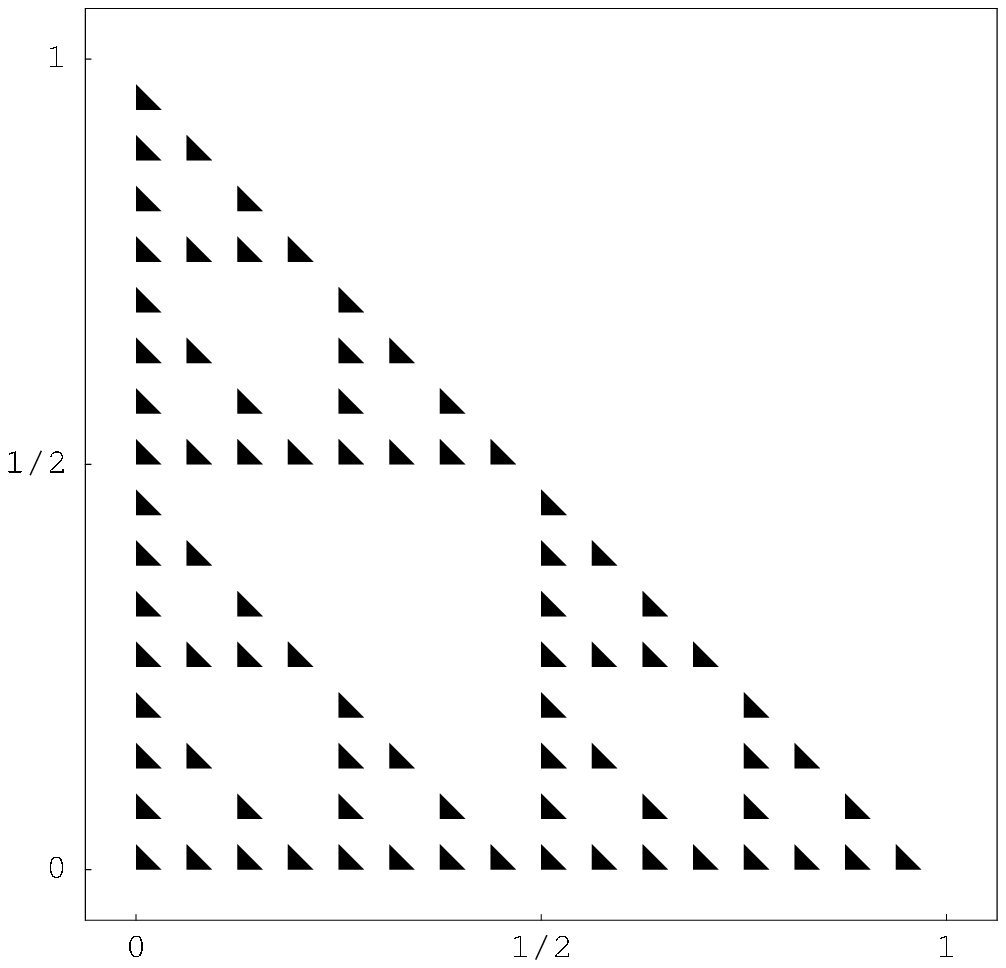}}
\put(2.125,1.625){\includegraphics[bb=88 4 376 292,width=\unitlength]{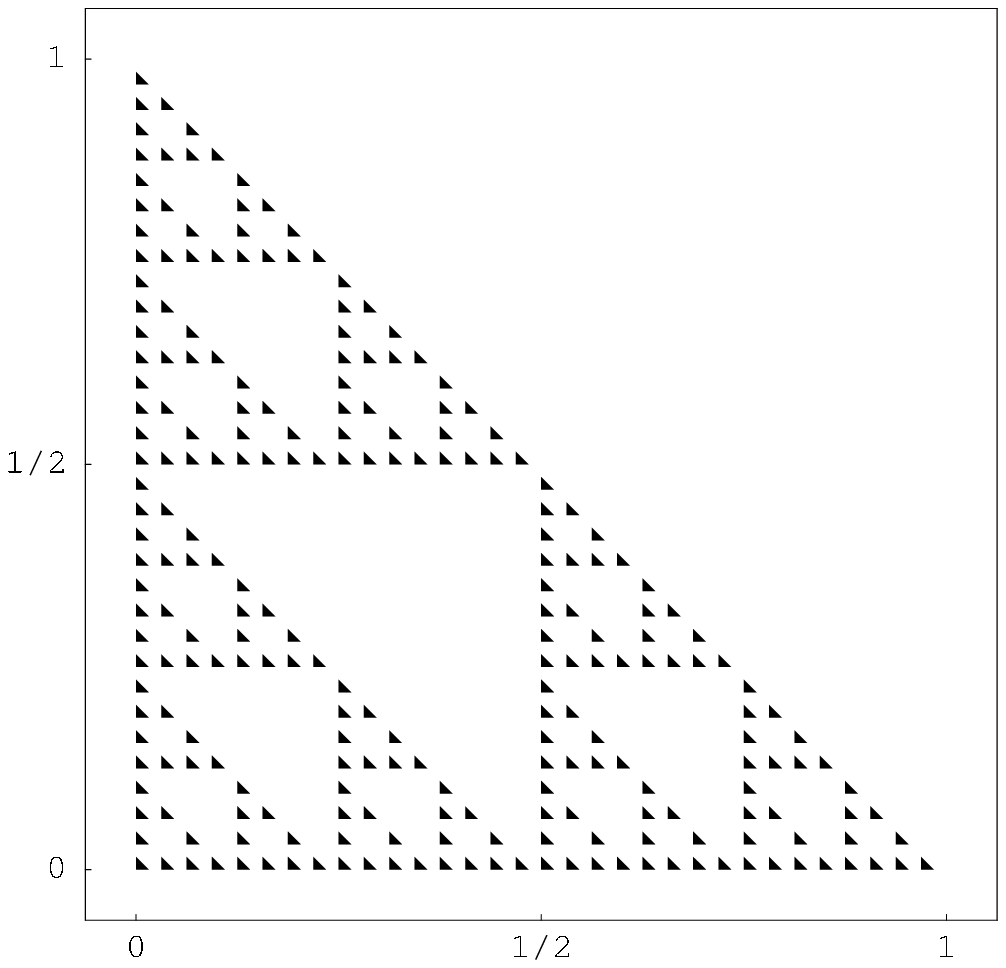}}
\put(0.60,1.625){\makebox(0,0)[t]{Fourth iteration}}
\put(1.6375,1.625){\makebox(0,0)[t]{Fifth iteration}}
\put(2.675,1.625){\makebox(0,0)[t]{Sixth iteration}}
\put(0,0){\includegraphics[bb=88 4 376 292,width=1.5\unitlength]{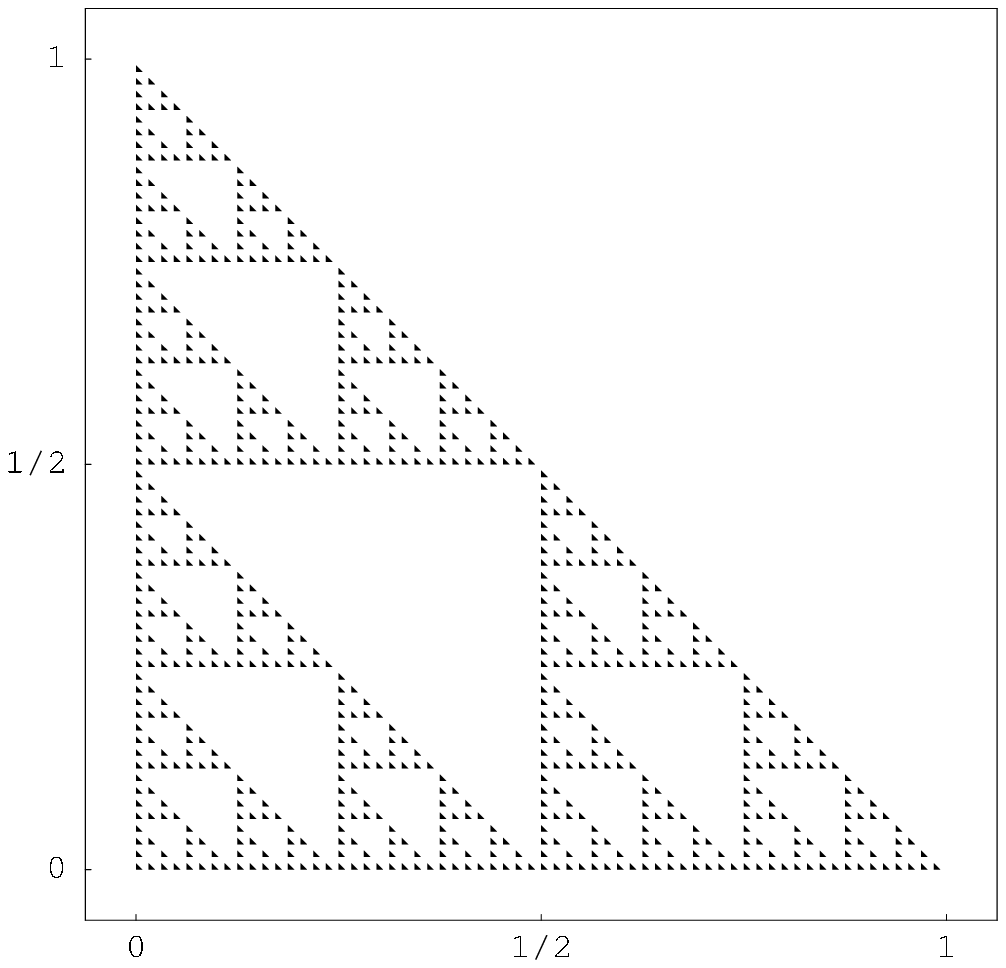}}
\put(1.625,0){\includegraphics[bb=88 4 376 292,width=1.5\unitlength]{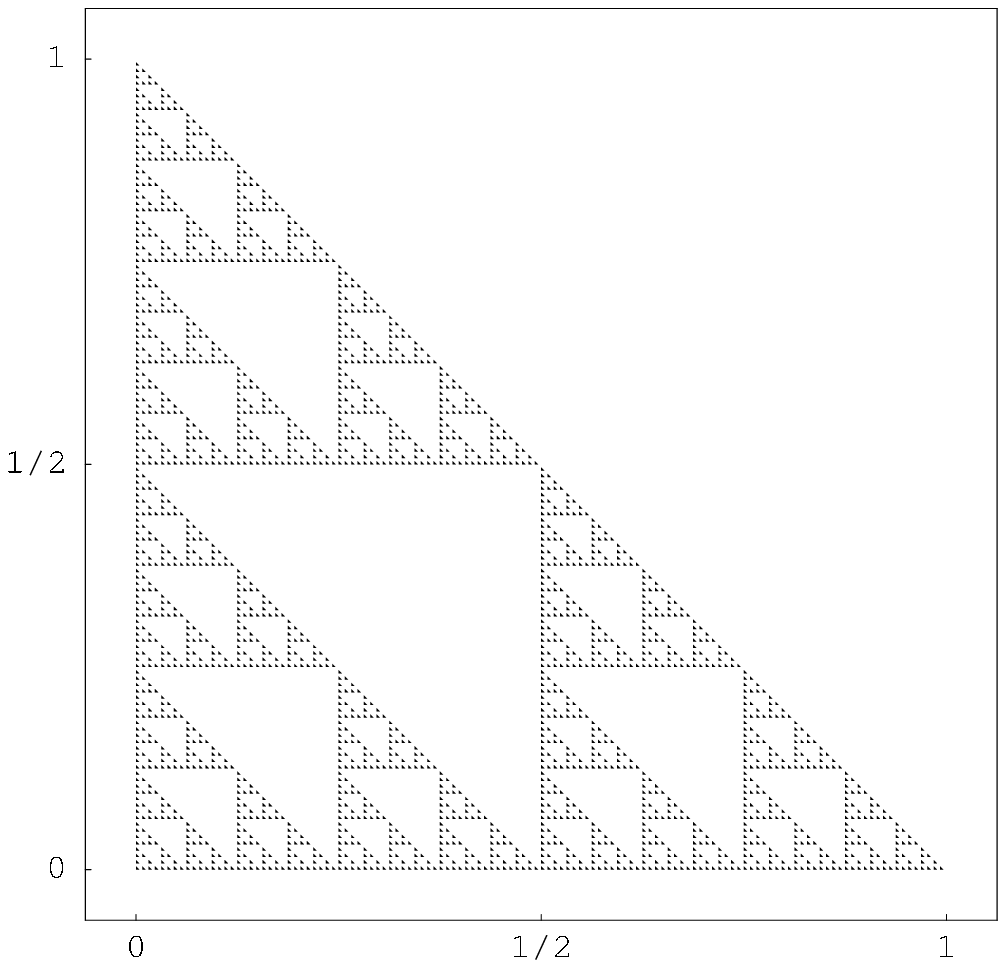}}
\put(0.825,0){\makebox(0,0)[t]{Seventh iteration}}
\put(2.45,0){\makebox(0,0)[t]{Eighth iteration}}
\end{picture}
\caption{The attractor $X_B$ for $R = \protect\rtwomatrix$, $B=\{(0,0),(1,0),(0,1)\}$}
\label{FigSierpTwoXB}
\end{figure}

\begin{figure}
\setlength{\unitlength}{0.48\textwidth}
\begin{picture}(2.083,2.722)(0,-0.125)
\put(0,2.0415){\includegraphics[bb=88 4 376 164,width=\unitlength]{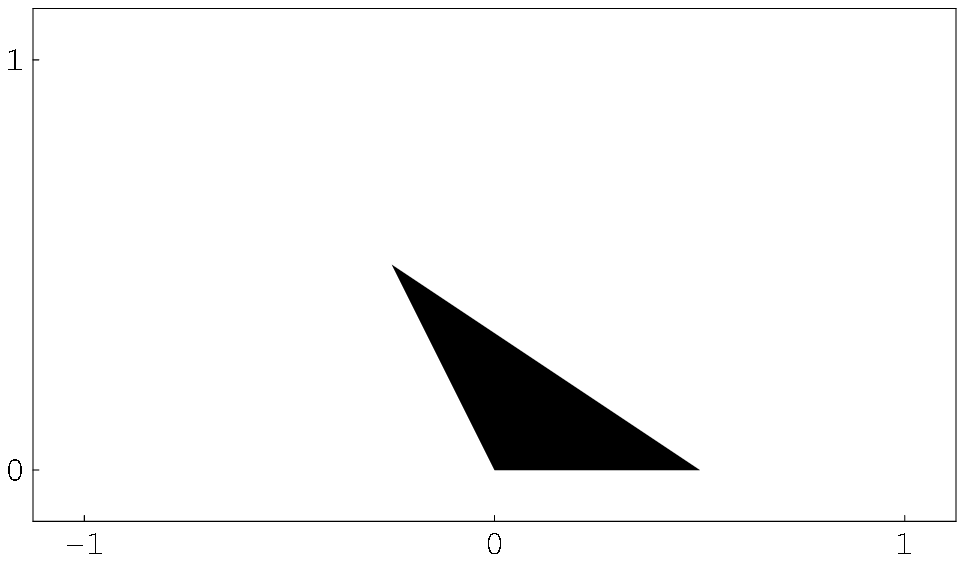}}
\put(1.083,2.0415){\includegraphics[bb=88 4 376 164,width=\unitlength]{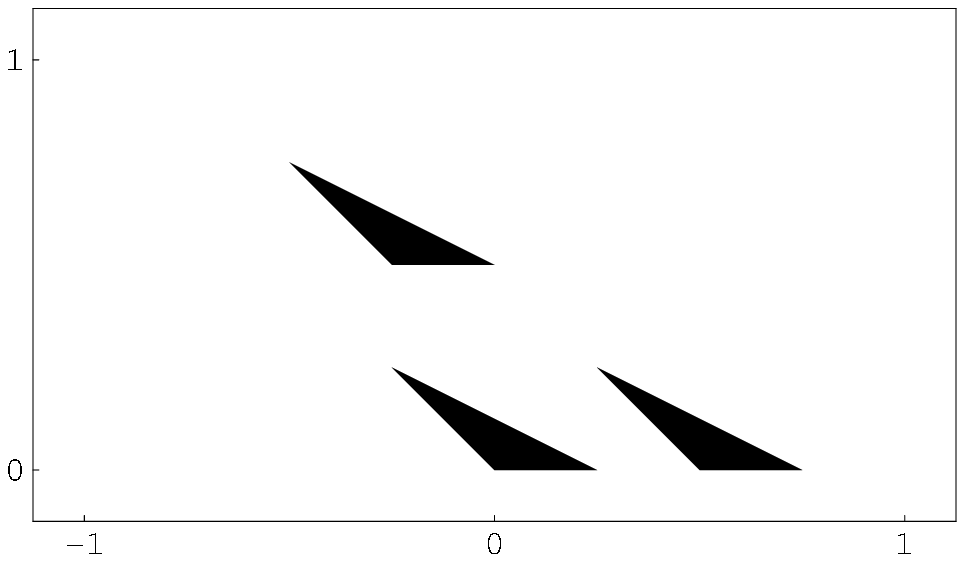}}
\put(0.55,2.0415){\makebox(0,0)[t]{First iteration}}
\put(1.633,2.0415){\makebox(0,0)[t]{Second iteration}}
\put(0,1.3610){\includegraphics[bb=88 4 376 164,width=\unitlength]{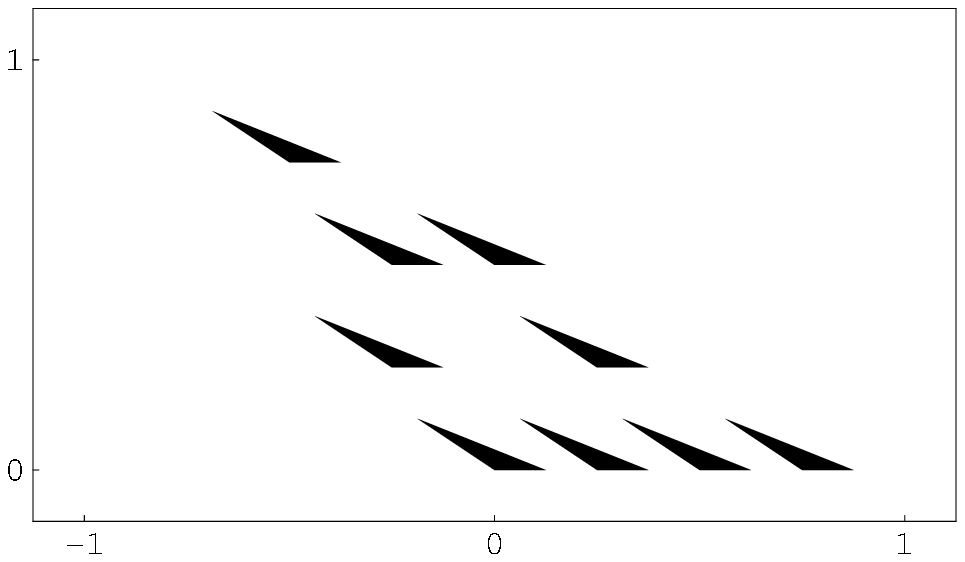}}
\put(1.083,1.3610){\includegraphics[bb=88 4 376 164,width=\unitlength]{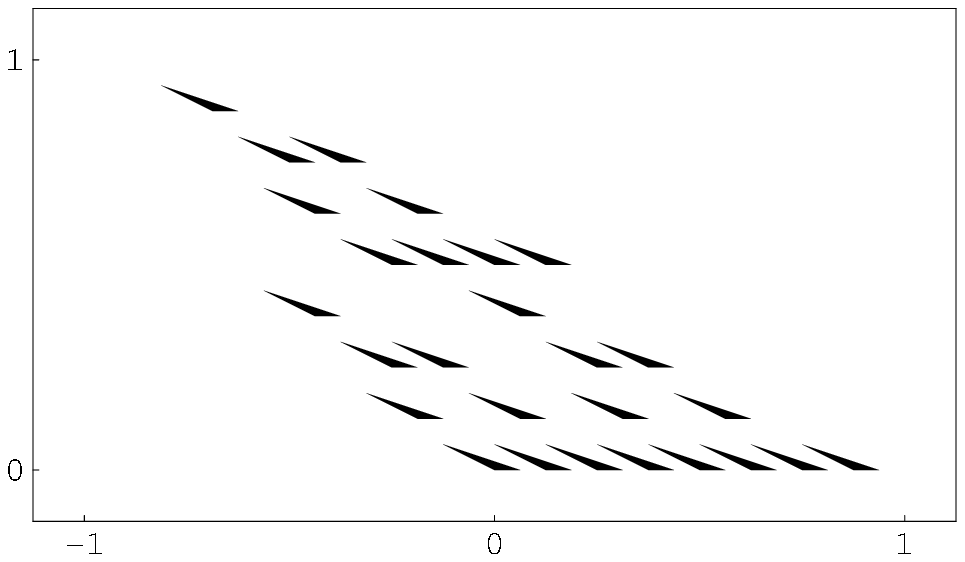}}
\put(0.55,1.3610){\makebox(0,0)[t]{Third iteration}}
\put(1.633,1.3610){\makebox(0,0)[t]{Fourth iteration}}
\put(0,0.6805){\includegraphics[bb=88 4 376 164,width=\unitlength]{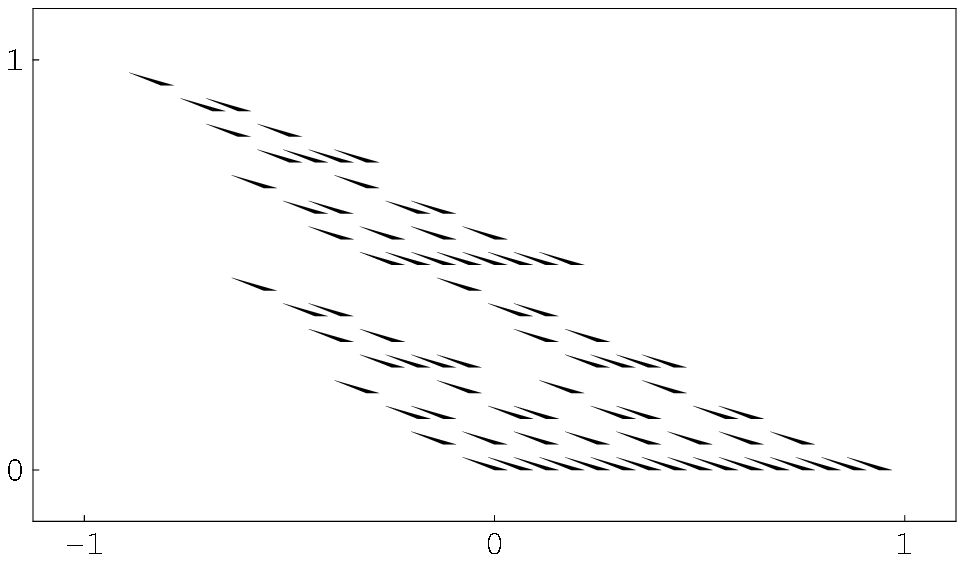}}
\put(1.083,0.6805){\includegraphics[bb=88 4 376 164,width=\unitlength]{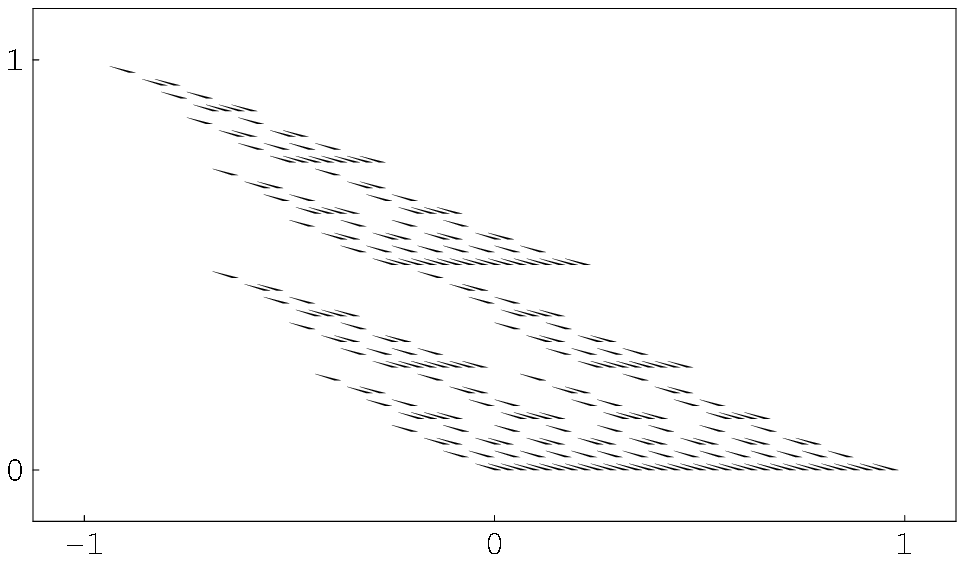}}
\put(0.55,0.6805){\makebox(0,0)[t]{Fifth iteration}}
\put(1.633,0.6805){\makebox(0,0)[t]{Sixth iteration}}
\put(0,0){\includegraphics[bb=88 4 376 164,width=\unitlength]{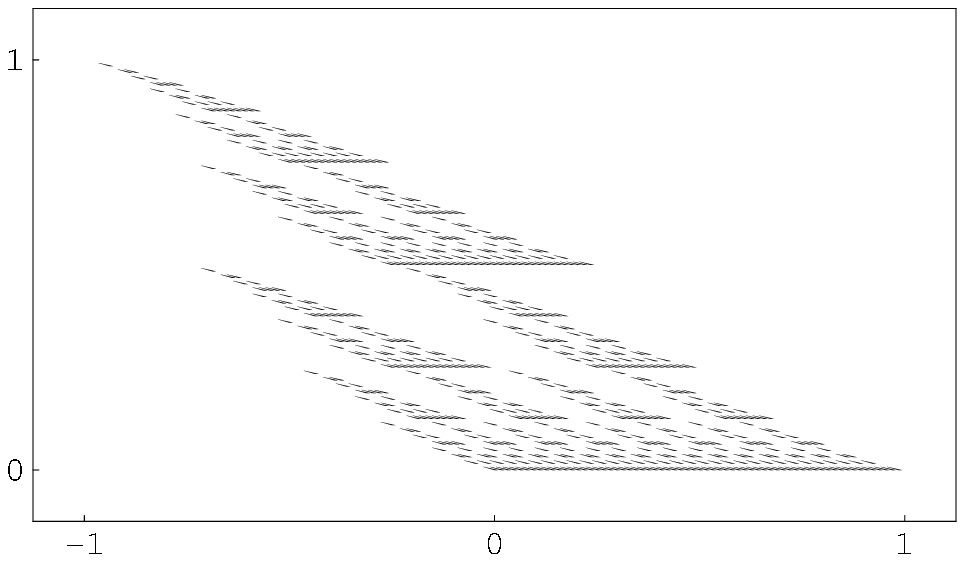}}
\put(1.083,0){\includegraphics[bb=88 4 376 164,width=\unitlength]{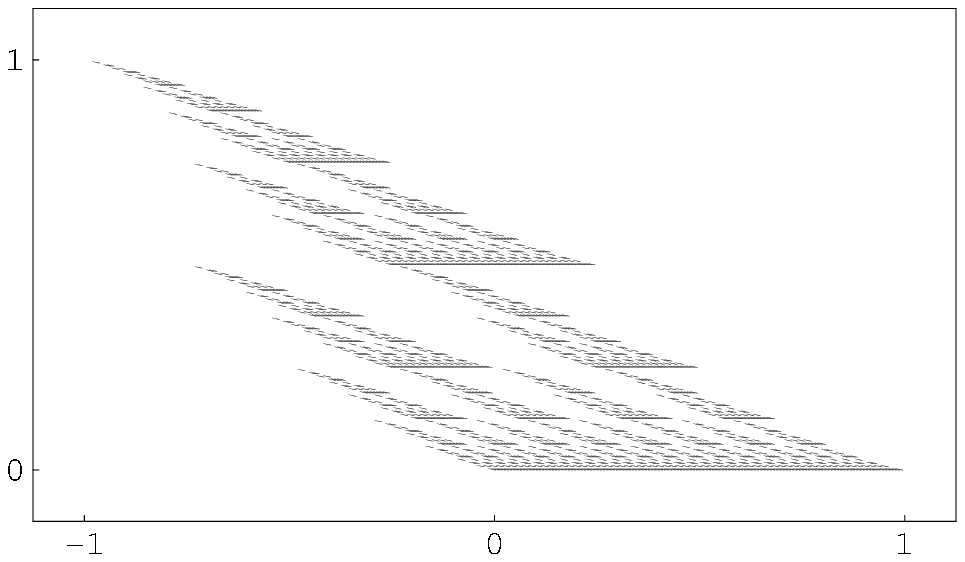}}
\put(0.55,0){\makebox(0,0)[t]{Seventh iteration}}
\put(1.633,0){\makebox(0,0)[t]{Eighth iteration}}
\end{picture}
\caption{The attractor $X_B$ for $R = \protect\rtwoonematrix$, $B=\{(0,0),(1,0),(0,1)\}$}
\label{FigSierpTwoOneXB}
\end{figure}

\begin{figure}
\setlength{\unitlength}{0.32\textwidth}
\begin{picture}(3.125,3.875)(0,-0.125)
\put(0.05,2.75){\includegraphics[bb=88 4 376 292,width=\unitlength]{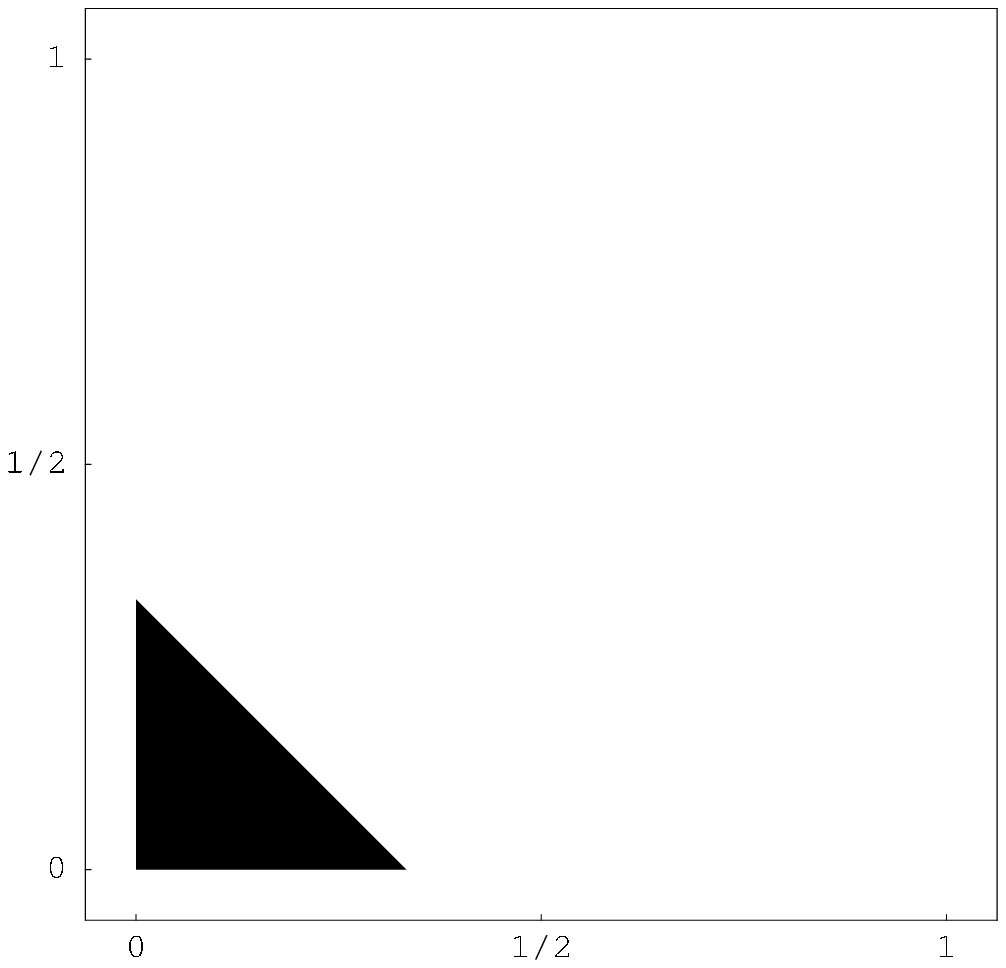}}
\put(1.0875,2.75){\includegraphics[bb=88 4 376 292,width=\unitlength]{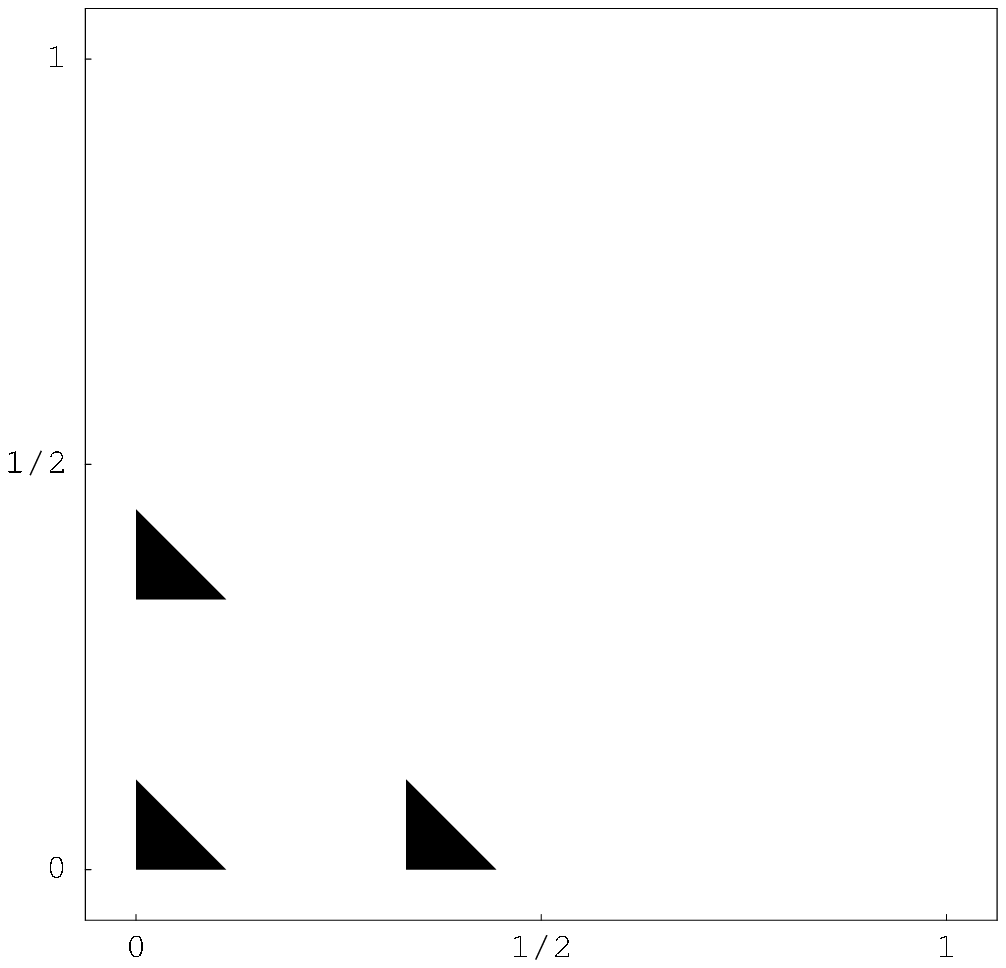}}
\put(2.125,2.75){\includegraphics[bb=88 4 376 292,width=\unitlength]{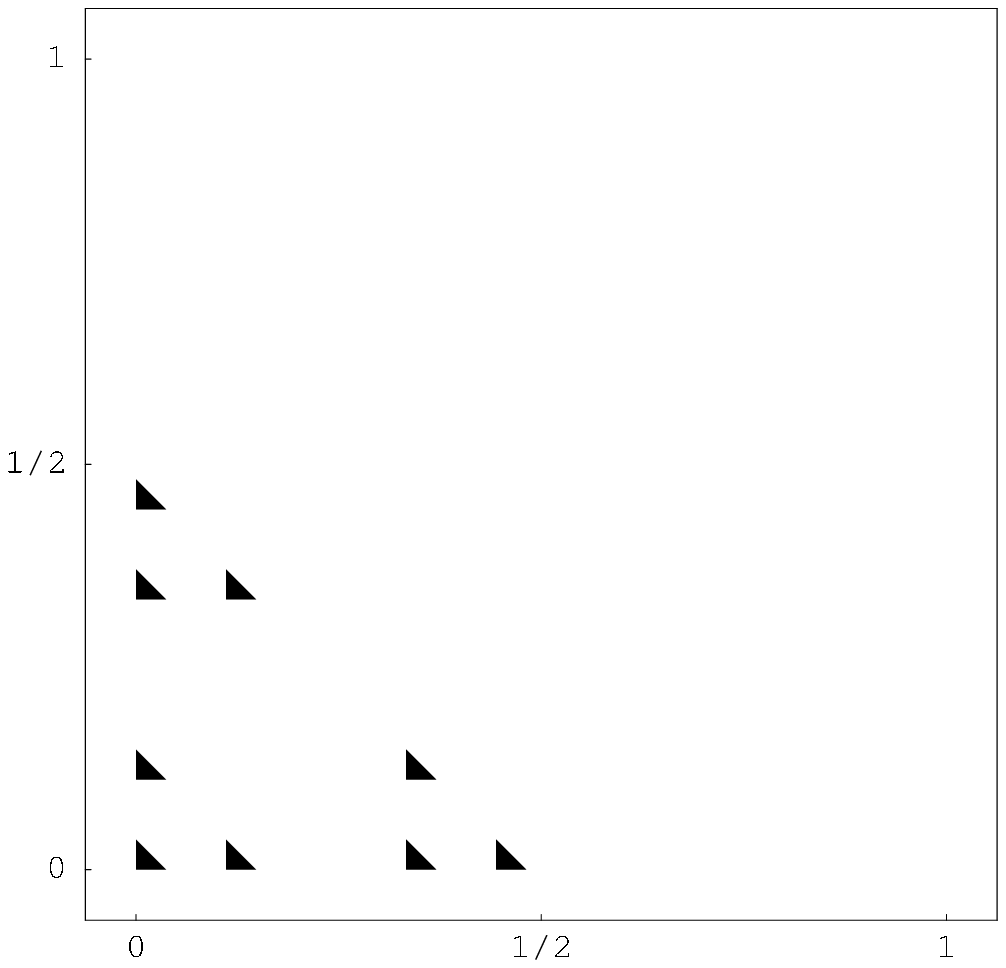}}
\put(0.60,2.75){\makebox(0,0)[t]{First iteration}}
\put(1.6375,2.75){\makebox(0,0)[t]{Second iteration}}
\put(2.675,2.75){\makebox(0,0)[t]{Third iteration}}
\put(0.05,1.625){\includegraphics[bb=88 4 376 292,width=\unitlength]{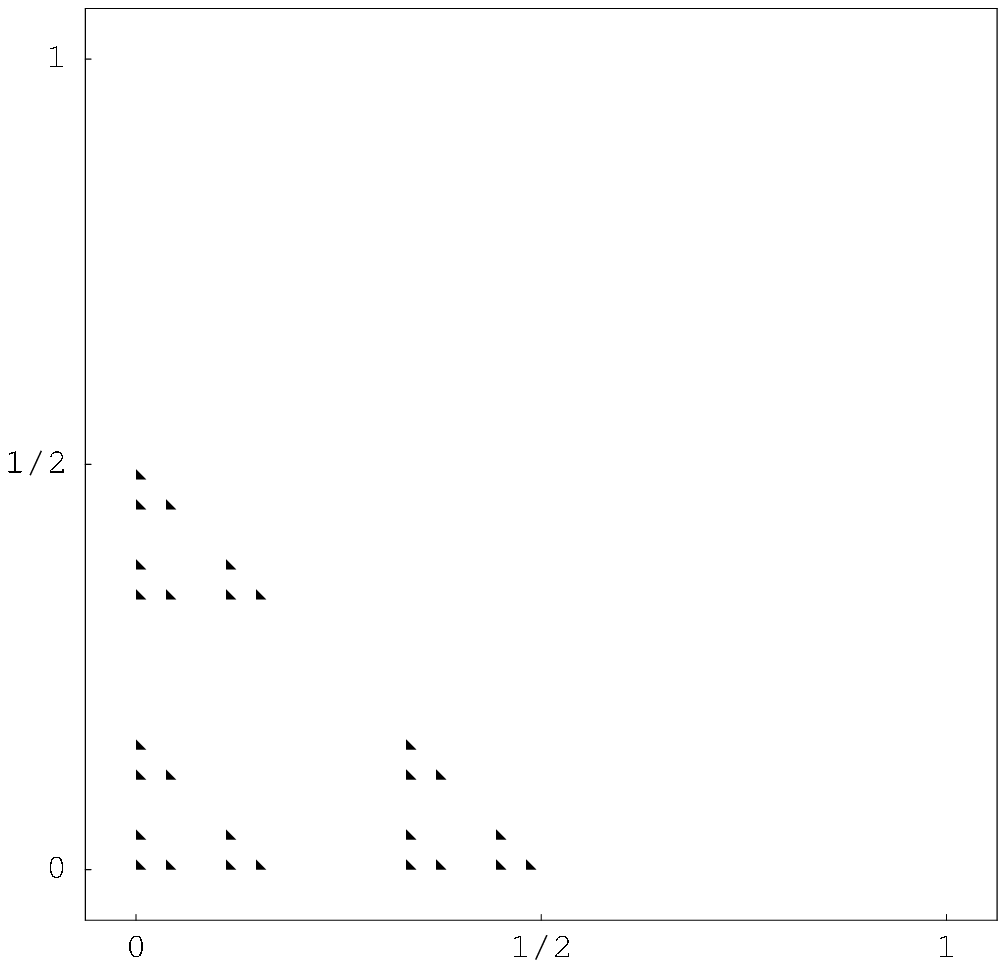}}
\put(1.0875,1.625){\includegraphics[bb=88 4 376 292,width=\unitlength]{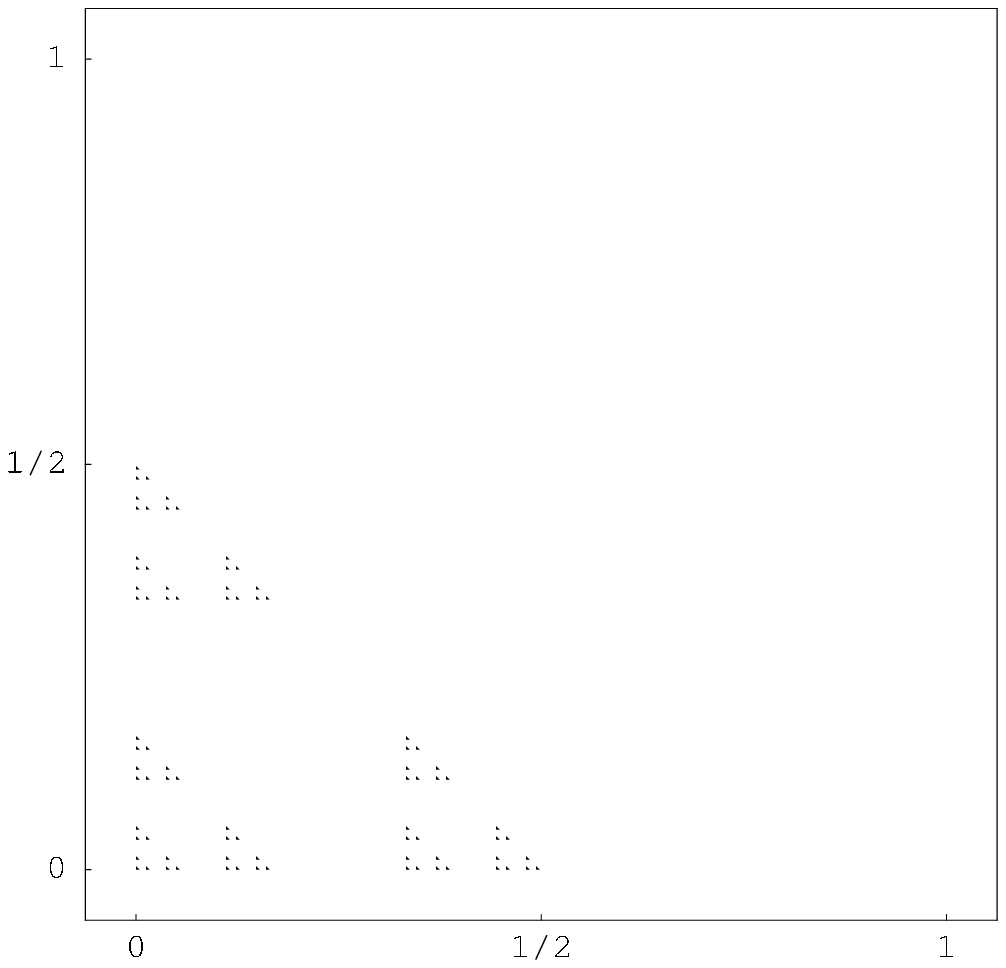}}
\put(2.125,1.625){\includegraphics[bb=88 4 376 292,width=\unitlength]{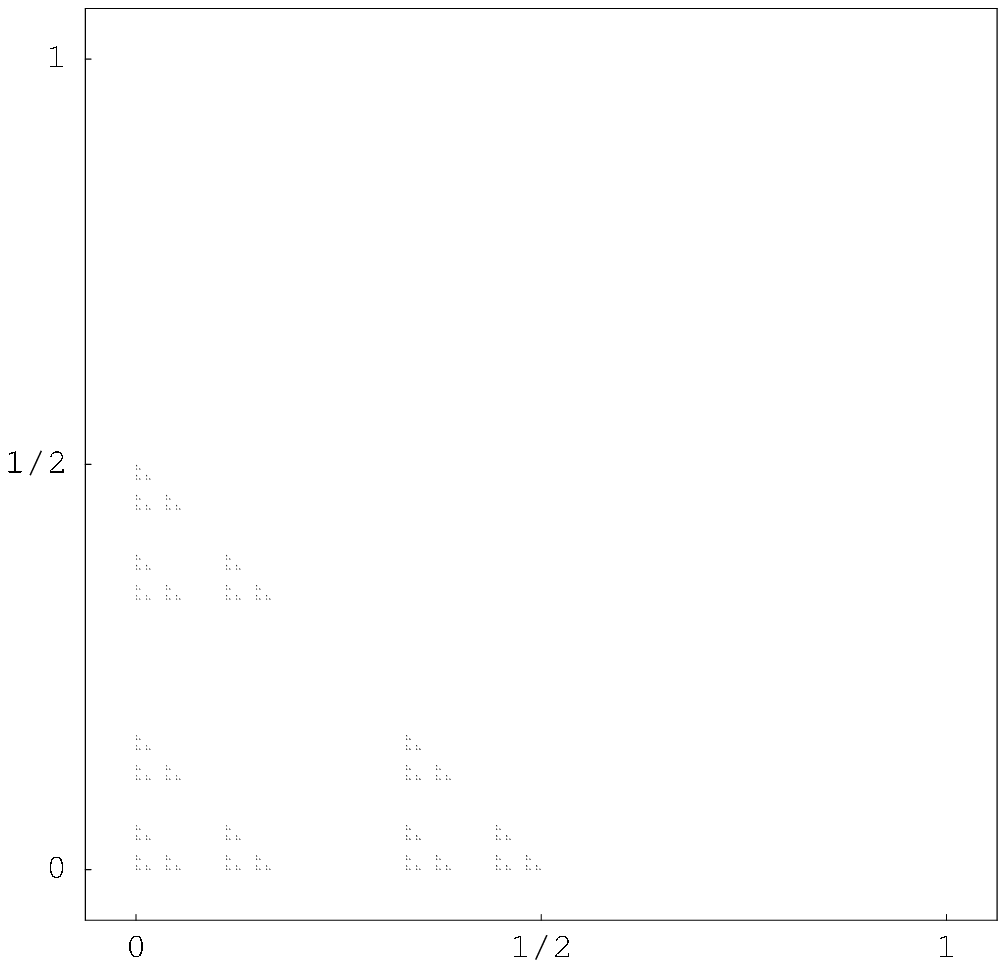}}
\put(0.60,1.625){\makebox(0,0)[t]{Fourth iteration}}
\put(1.6375,1.625){\makebox(0,0)[t]{Fifth iteration}}
\put(2.675,1.625){\makebox(0,0)[t]{Sixth iteration}}
\put(0,0){\includegraphics[bb=88 4 376 292,width=1.5\unitlength]{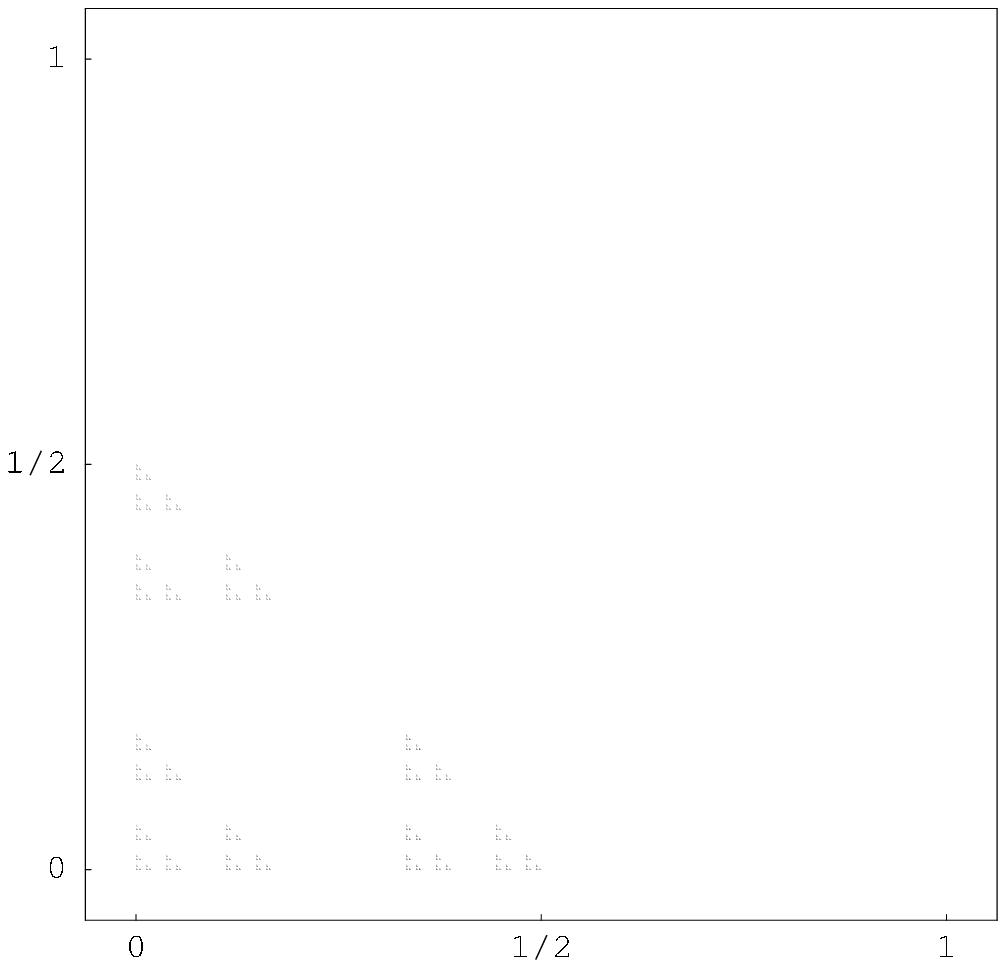}}
\put(1.625,0){\includegraphics[bb=88 4 376 292,width=1.5\unitlength]{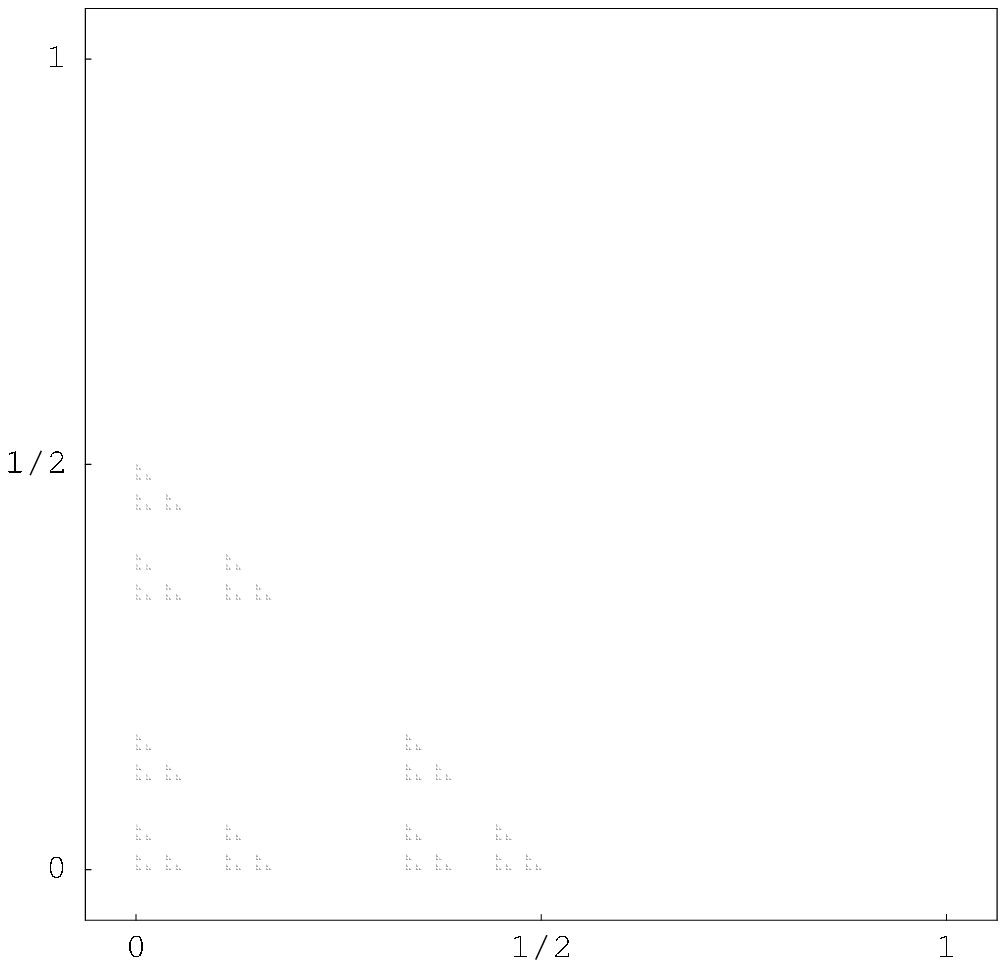}}
\put(0.825,0){\makebox(0,0)[t]{Seventh iteration}}
\put(2.45,0){\makebox(0,0)[t]{Eighth iteration}}
\end{picture}
\caption{The attractor $X_B$ for $R = \protect\rthreematrix$, $B=\{(0,0),(1,0),(0,1)\}$}
\label{FigSierpThreeXB}
\end{figure}

     In Figures \ref{FigSierpTwoXB}, \ref{FigSierpTwoOneXB}, and \ref{FigSierpThreeXB}
we have included planar Sierpinski examples
corresponding to a fixed configuration of the three vectors 
$B=\{(0,0),(1,0),(0,1)\}$, and scaling matrices 
$R = \rtwomatrix$, $R = \rtwoonematrix$, and $R = \rthreematrix$, respectively. 

The three cases serve to illustrate that spectral properties of the measure
$\mu_B$, or rather $\mu_{R,B}$, depend on the choice of scaling matrix in an
essential way: For the first two examples $\mu_B$ in fact is not a spectral
measure, while it is in the last example.


\begin{figure}
\setlength{\unitlength}{0.32\textwidth}
\begin{picture}(3.125,3.875)(0,-0.125)
\put(0.05,2.75){\includegraphics[bb=79 4 367 292,width=\unitlength]{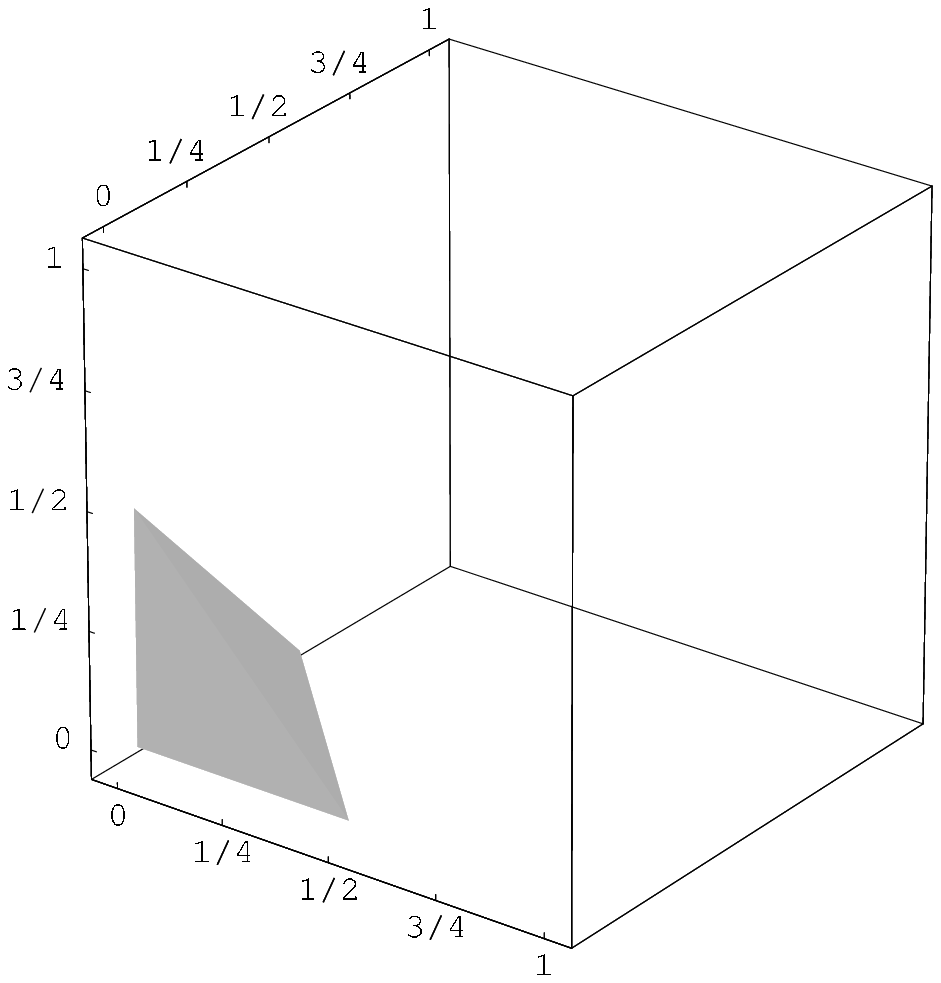}}
\put(1.0875,2.75){\includegraphics[bb=79 4 367 292,width=\unitlength]{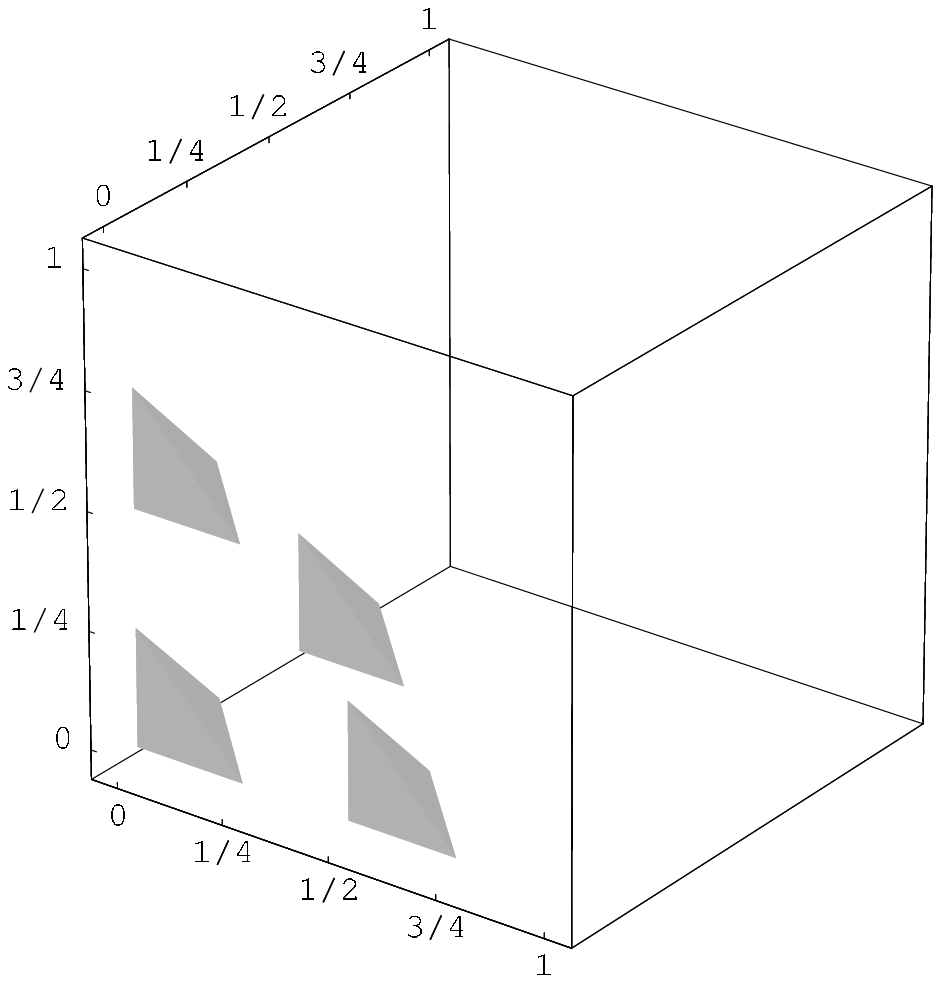}}
\put(2.125,2.75){\includegraphics[bb=79 4 367 292,width=\unitlength]{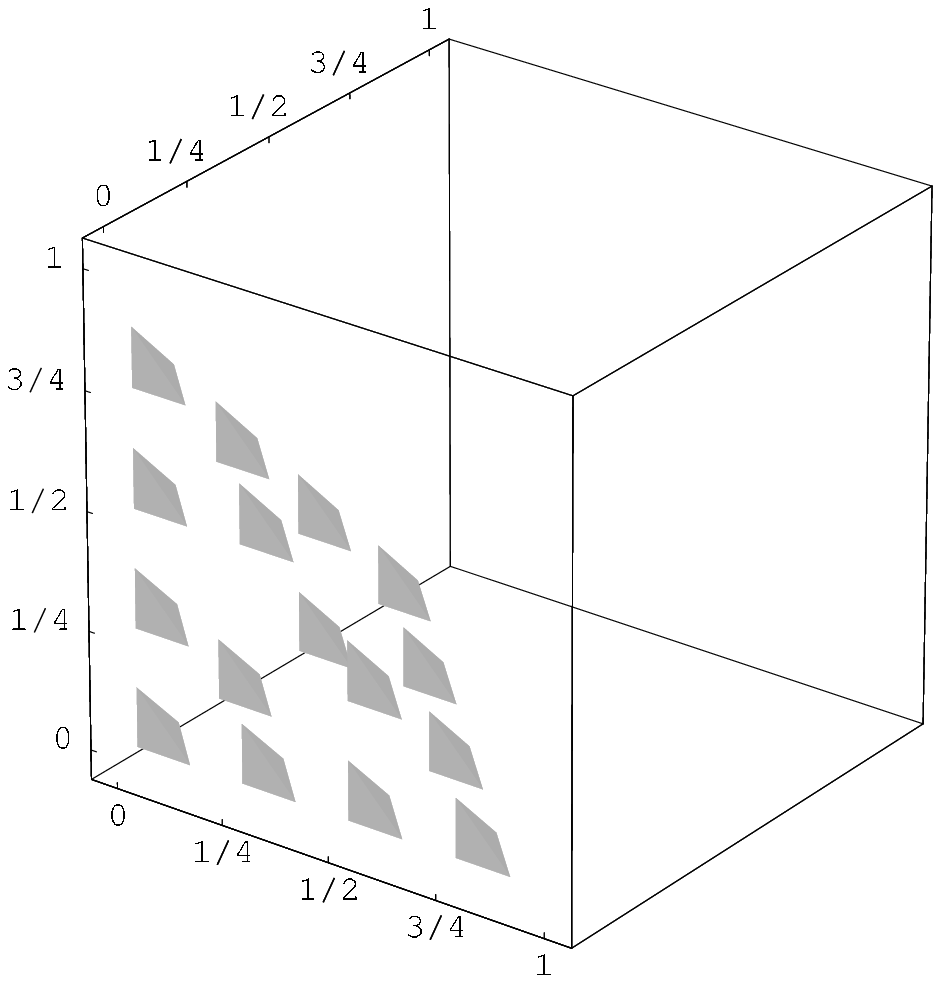}}
\put(0.60,2.75){\makebox(0,0)[t]{First iteration}}
\put(1.6375,2.75){\makebox(0,0)[t]{Second iteration}}
\put(2.675,2.75){\makebox(0,0)[t]{Third iteration}}
\put(0.05,1.625){\includegraphics[bb=79 4 367 292,width=\unitlength]{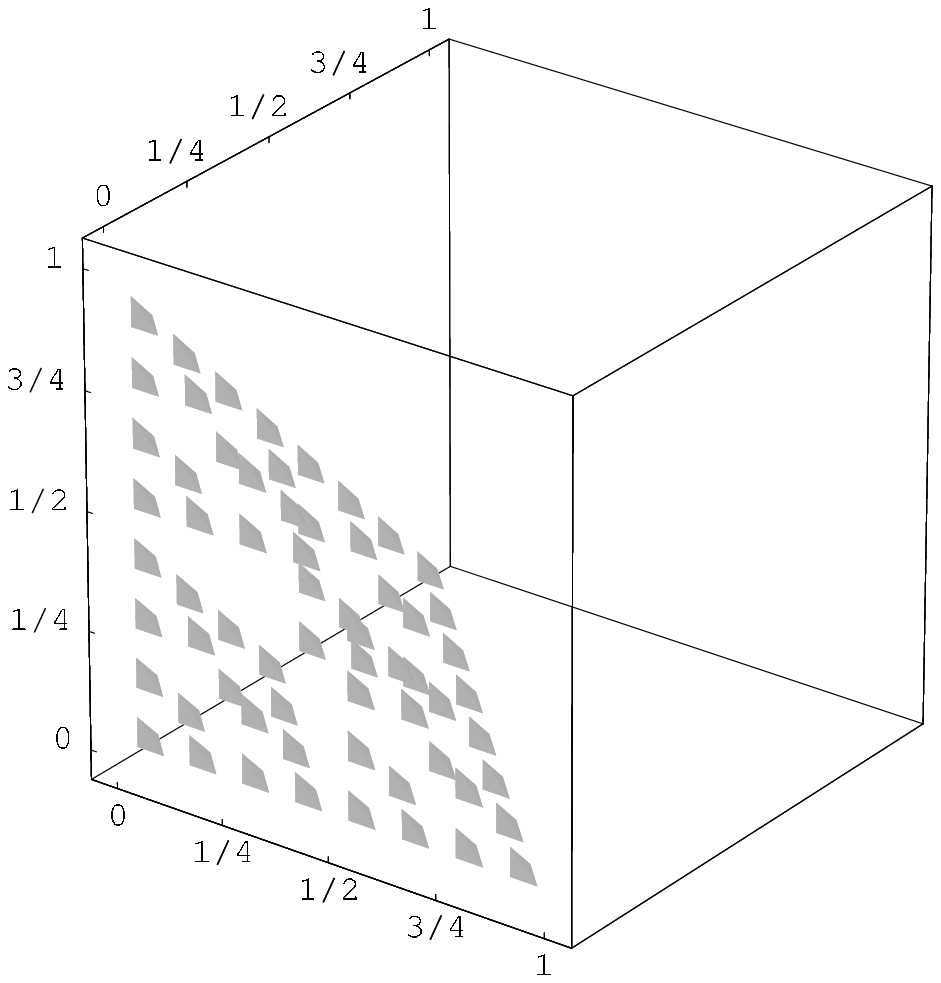}}
\put(1.0875,1.625){\includegraphics[bb=79 4 367 292,width=\unitlength]{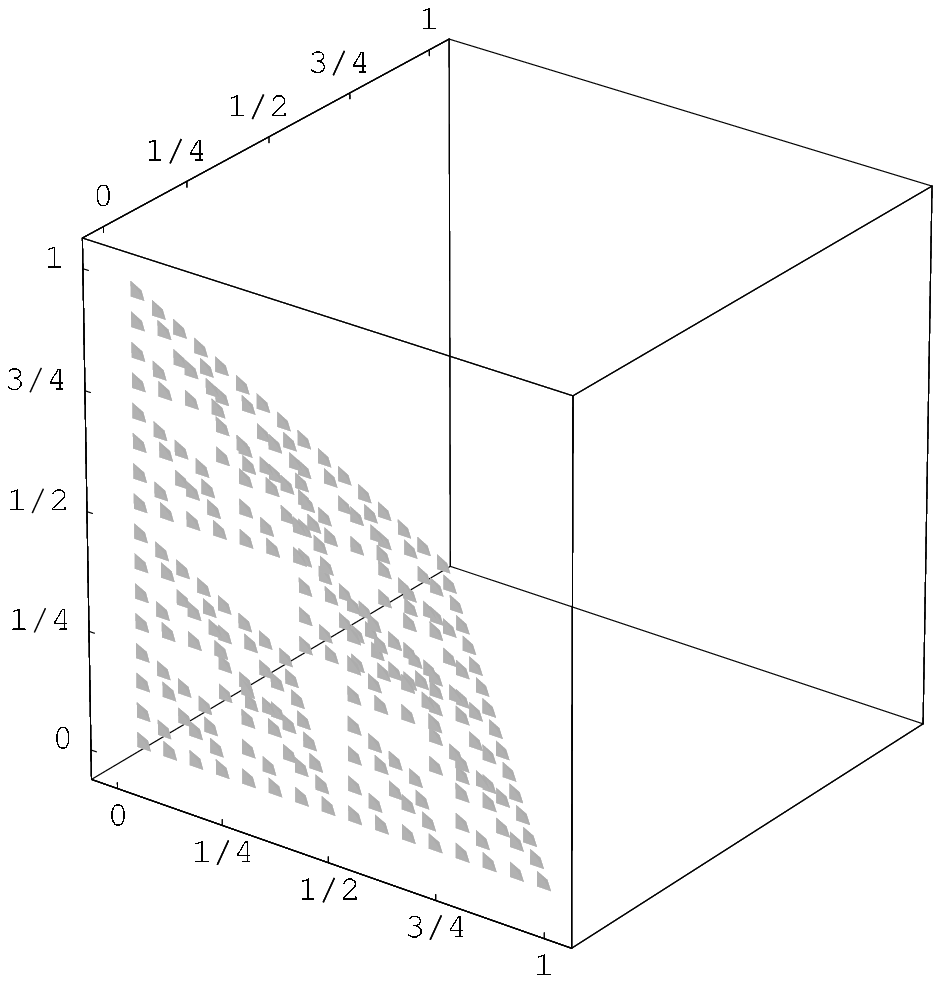}}
\put(2.125,1.625){\includegraphics[bb=79 4 367 292,width=\unitlength]{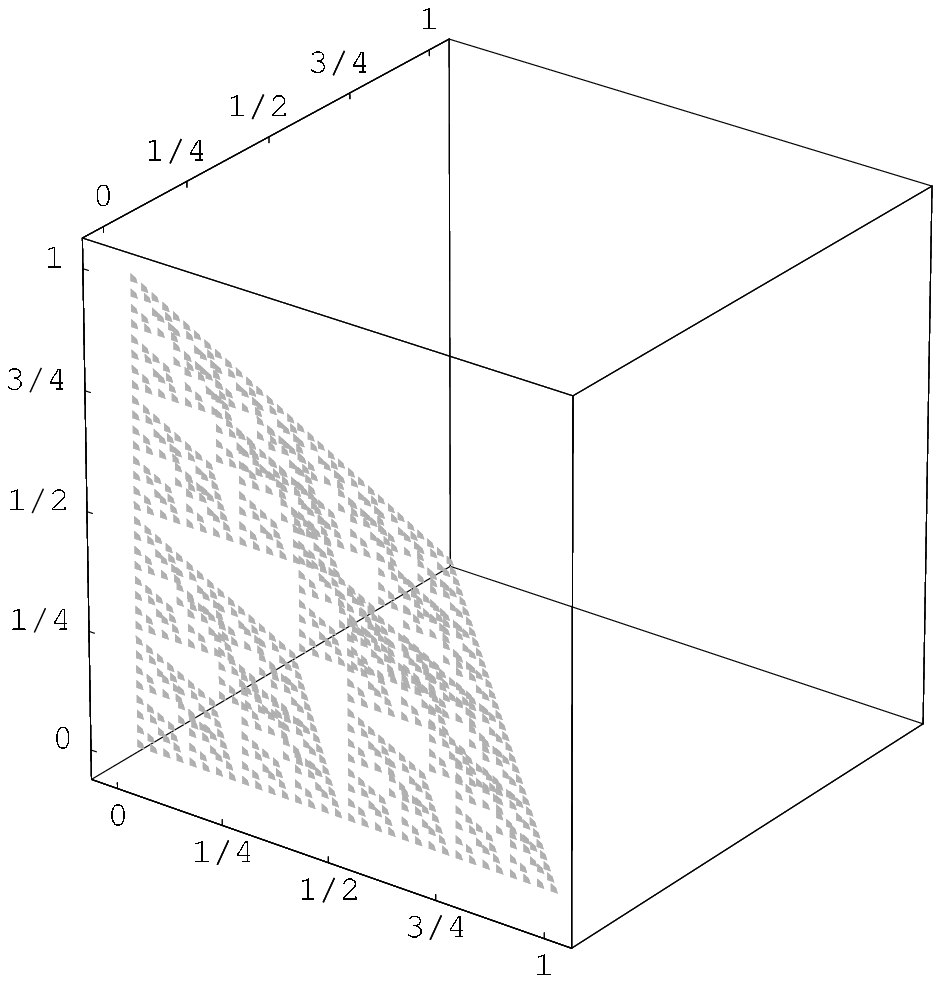}}
\put(0.60,1.625){\makebox(0,0)[t]{Fourth iteration}}
\put(1.6375,1.625){\makebox(0,0)[t]{Fifth iteration}}
\put(2.675,1.625){\makebox(0,0)[t]{Sixth iteration}}
\put(0.825,0){\makebox(0,0)[t]{Seventh iteration}}
\put(2.45,0){\makebox(0,0)[t]{Eighth iteration}}
\end{picture}
\caption{The attractor $X_B$ for $R = \protect\rtwomatrixthreed$,
$B=\left\{\protect\threevectoro,\protect\threevectorx,\protect\threevectory,\protect\threevectorz\right\}$}
\label{FigSierp3DXB}
\end{figure}

\begin{figure}
\setlength{\unitlength}{0.32\textwidth}
\begin{picture}(3.125,3.875)(0,-0.125)
\put(0.05,2.75){\includegraphics[bb=79 4 367 292,width=\unitlength]{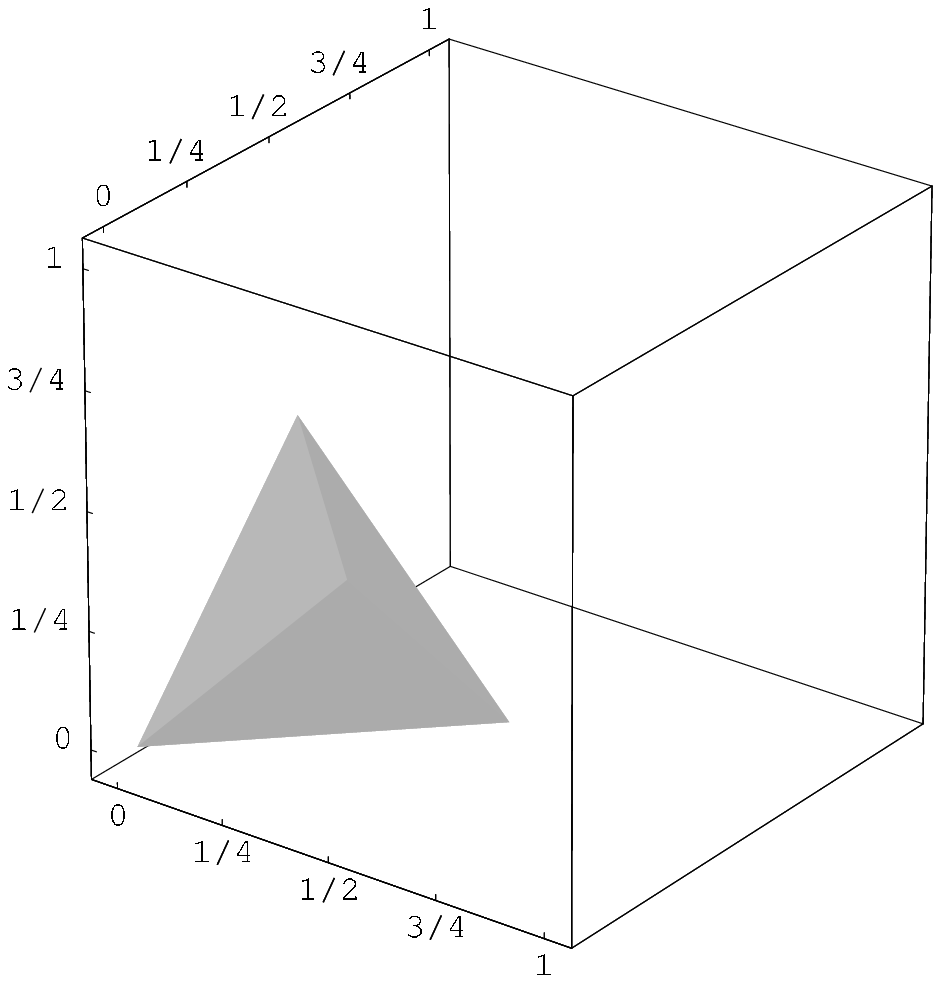}}
\put(1.0875,2.75){\includegraphics[bb=79 4 367 292,width=\unitlength]{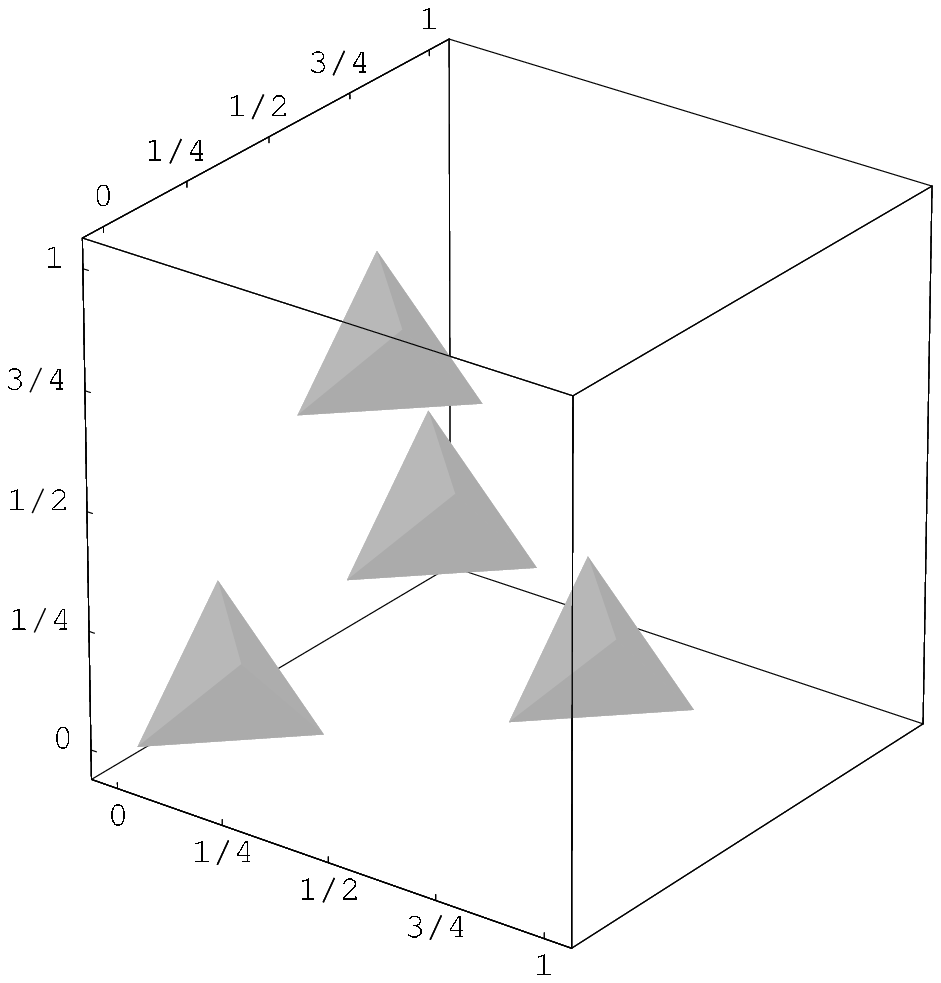}}
\put(2.125,2.75){\includegraphics[bb=79 4 367 292,width=\unitlength]{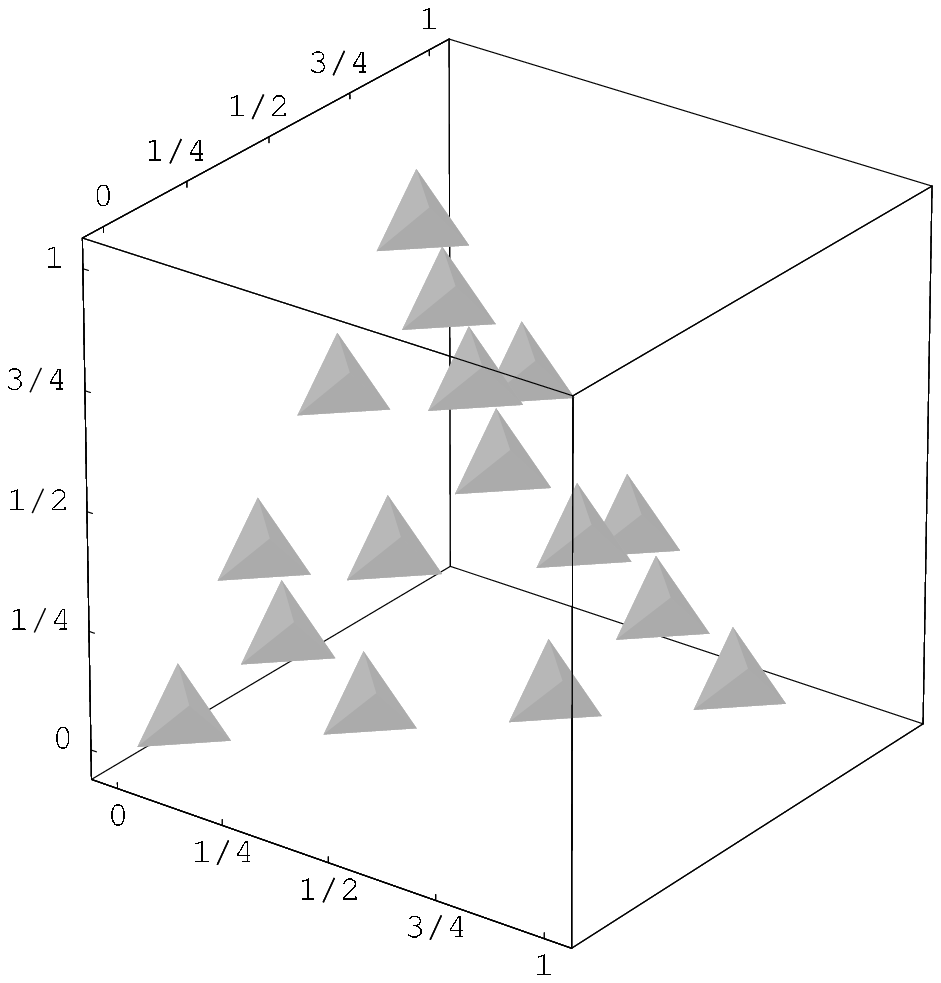}}
\put(0.60,2.75){\makebox(0,0)[t]{First iteration}}
\put(1.6375,2.75){\makebox(0,0)[t]{Second iteration}}
\put(2.675,2.75){\makebox(0,0)[t]{Third iteration}}
\put(0.05,1.625){\includegraphics[bb=79 4 367 292,width=\unitlength]{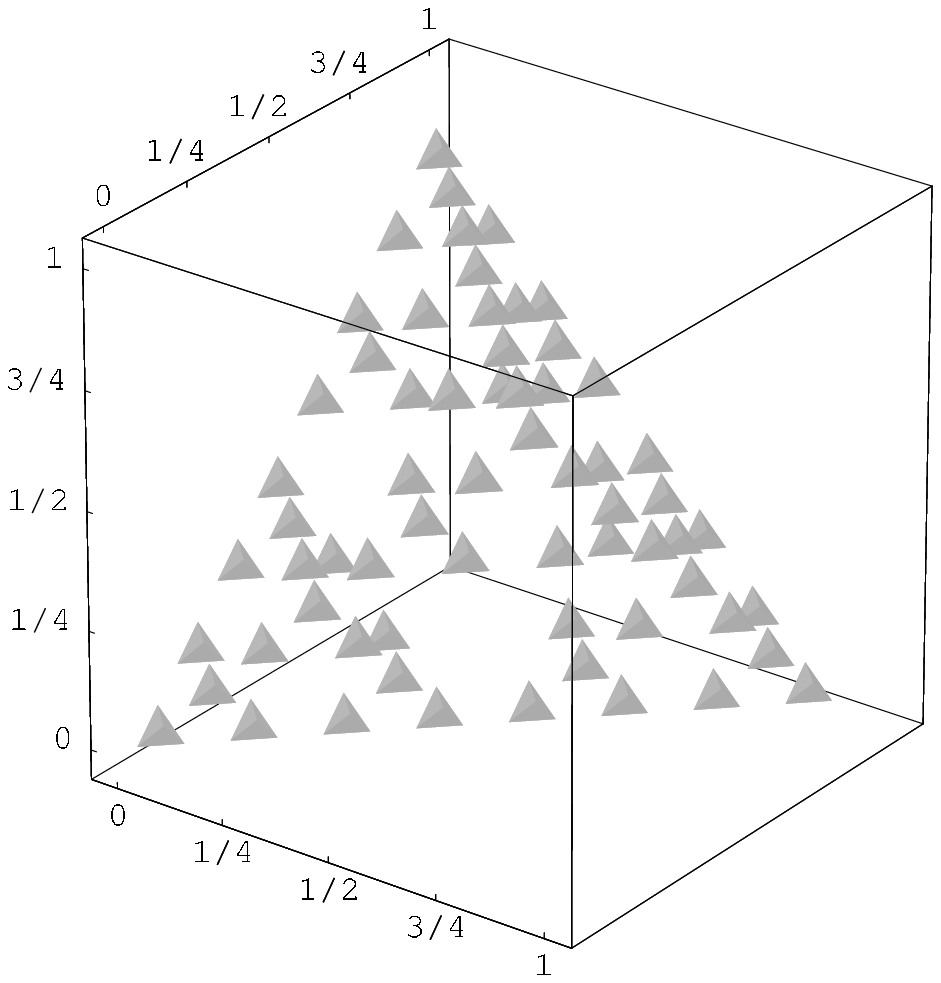}}
\put(1.0875,1.625){\includegraphics[bb=79 4 367 292,width=\unitlength]{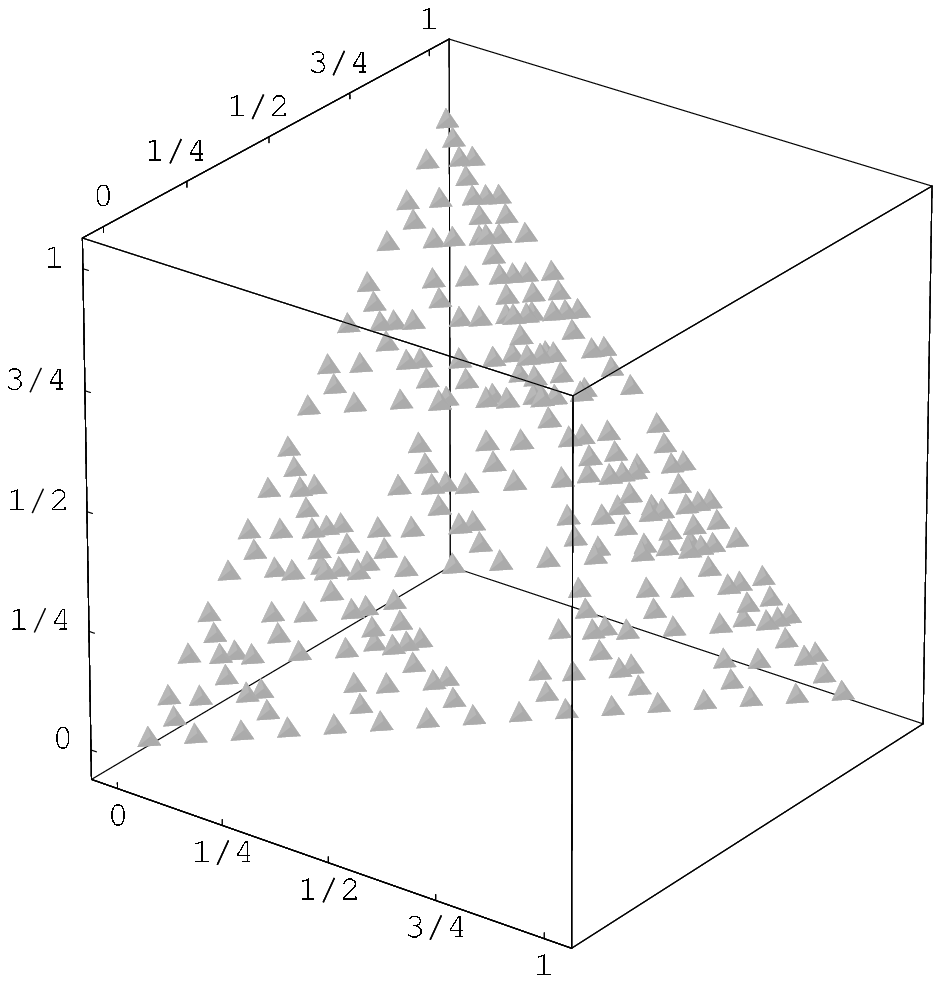}}
\put(2.125,1.625){\includegraphics[bb=79 4 367 292,width=\unitlength]{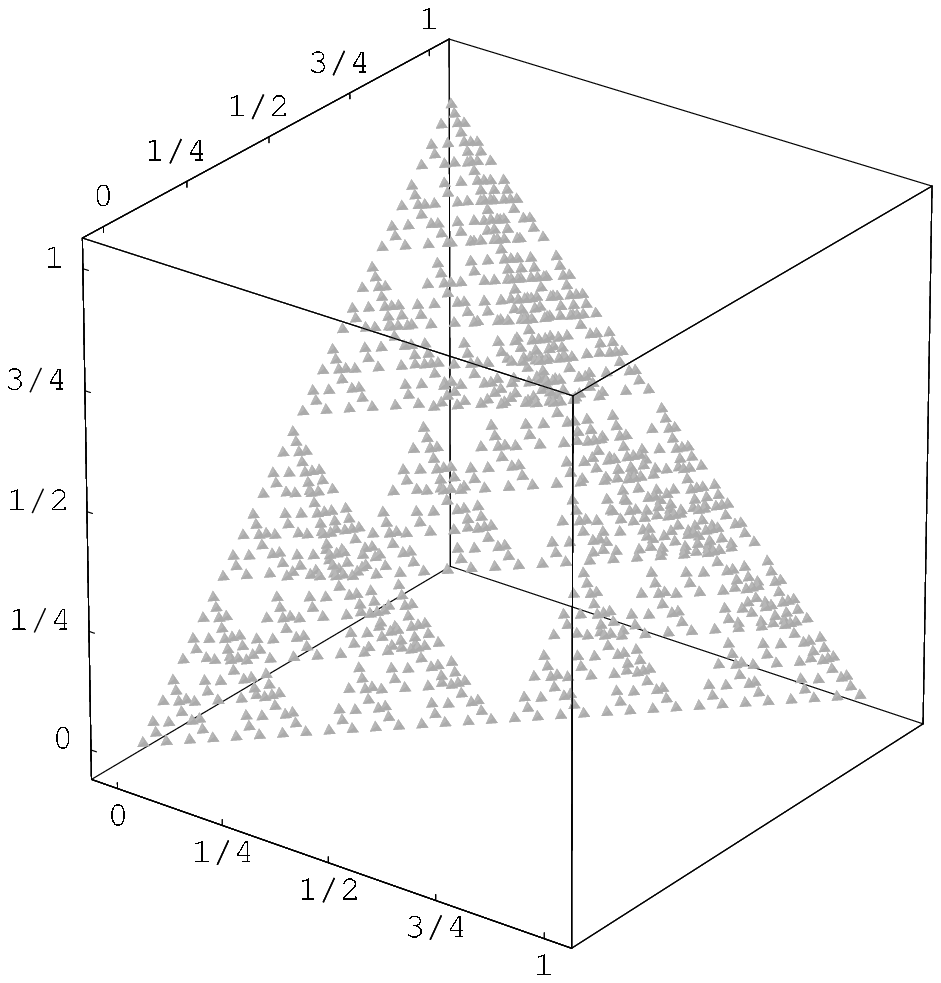}}
\put(0.60,1.625){\makebox(0,0)[t]{Fourth iteration}}
\put(1.6375,1.625){\makebox(0,0)[t]{Fifth iteration}}
\put(2.675,1.625){\makebox(0,0)[t]{Sixth iteration}}
\put(0.825,0){\makebox(0,0)[t]{Seventh iteration}}
\put(2.45,0){\makebox(0,0)[t]{Eighth iteration}}
\end{picture}
\caption{Dual attractor $X_L$ for $S = \protect\rtwomatrixthreed$,
$L=\left\{\protect\threevectoro,\protect\threevectorxy,\protect\threevectorxz,\protect\threevectoryz\right\}$}
\label{FigSierp3DXL}
\end{figure}

In Figures \ref{FigSierp3DXB} and \ref{FigSierp3DXL} below we sketch the dual system of attractors $X_B$ and $X_L$ for the configuration of $(R, B, L)$ in the case of
part (\ref{thsier(3)}) of Theorem \ref{thsier} when $p$ is even (specifically, $p=2$).
For this particular
Sierpinski configuration, we arrive at a 3D example where $\mu_B$ is a
spectral measure.

\par
We consider the following type of affine iterated function systems: $R=pI_d$, with $p\in\bz$, $p\geq 2$, and $B$ consists of $0$ and $d$ linearly independent vectors in $\br^d$. For example, when $d=2$ $R=2I_2$ and $B=\{(0,0),(1,0),(1/2,\sqrt{3}/2)\}$, we get the classical Sierpinski gasket.
\par
We are interested for which of these examples $\mu_B$ is a spectral measure. 
\par
Note that, using a change of variable as in Lemma \ref{lemtran}, we can always take $B=\{0,e_1,\dots,e_d\}$ where $e_k$ are the cannonical vectors $e_k=(0,\dots,0,1,0,\dots,0)$ with $1$ in the $k$-th position. Indeed, just take $V$ in Lemma \ref{lemtran} to be the matrix that maps the cannonical vectors into the vectors in $B$.  
\begin{theorem}\label{thsier}
Let $R=pI_d$, $B=\{0,e_1,\dots,e_d\}$, and let $\mu_B$ be the associated invariant measure.
\begin{enumerate}
\item\label{thsier(1)} If the dimension $d$ is $1$, then $\mu_B$ is a spectral measure if and only if $p$ is even. \par
If $p$ is odd, then there are no $3$ mutually orthogonal exponential functions in $L^2(\mu_B)$.\par
If $p$ is even then one can take $L=\{0,p/2\}$, to obtain a Hadamard triple $(R,B,L)$, and
\begin{itemize}
\item\label{thsier(1)(1)} If $p=2$, then there are two $W_B$-cycles, $\{0\}$ and $\{1\}$ and the spectrum of $\mu_B$ is $\bz$;
\item\label{thsier(1)(2)} If $p>2$, then there is only one $W_B$-cycle $\{0\}$ and the spectrum of $\mu_B$ is $$\{\frac p2\sum_{k=0}^np^ka_k\,|\,a_k\in\{0,1\},n\in\bn\}.$$
\end{itemize}
\item\label{thsier(2)} If the dimension $d$ is $2$, then $\mu_B$ is a spectral measure if and only if $p$ is a multiple of $3$. 
\par
If $p$ is not a multiple of $3$ then there are no $4$ mutually orthogonal exponential functions in $L^2(\mu_B)$. \par
If $p$, is a multiple of $3$ then one can take $L:=\{(0,0),(\frac{2p}3,-\frac{2p}3),(-\frac{2p}3,\frac{2p}3)\}$ to obtain a Hadamard triple $(R,B,L)$. There is only one $W_B$-cycle $\{(0,0)\}$ and the spectrum of $\mu_B$ is $$\{\frac{2p}3\sum_{k=0}^np^ka_k\,|\,a_k\in\{(0,0),(-1,1),(1,-1)\},n\in\bn\}.$$
\item\label{thsier(3)}
If the dimension $d$ is $3$, then $\mu_B$ is a spectral measure if and only if $p$ is even.\par
If $p$ is odd then there are at most $256$ mutually orthogonal exponential functions in $L^2(\mu_B)$.\par
If $p$ is even, then one can take $L:=\{(0,0,0),(\frac p2,\frac p2,0),(0,\frac p2,\frac p2),(\frac p2,0,\frac p2)\}$, to obtain a Hadamard triple $(R,B,L)$, and
\begin{itemize}
\item\label{thsier(3)(1)} If $p=2$, then there are four $W_B$-cycles, $\{(0,0,0),(1,1,0),(0,1,1),(1,0,1)\}$, and the spectrum of $\mu_B$ is the union of the following sets:
\begin{align*}
\Lambda(0,0,0)&=\{\sum_{k=0}^{n-1}2^kl_k\,|\,l_k\in L,n\in\bn\}\\
\Lambda(1,1,0)&=\{2^n(1,1,0)-\sum_{k=0}^{n-1}2^kl_k\,|\,l_k\in L,n\in\bn\}\\
\Lambda(1,0,1)&=\{2^n(1,0,1)-\sum_{k=0}^{n-1}2^kl_k\,|\,l_k\in L,n\in\bn\}\\
\Lambda(0,1,1)&=\{2^n(0,1,1)-\sum_{k=0}^{n-1}2^kl_k\,|\,l_k\in L,n\in\bn\}.
\end{align*}
\item\label{thsier(3)(2)} If $p>2$ and is even, then there is only one $W_B$-cycle $\{0,0,0\}$ and the spectrum of $\mu_B$ is 
$$\Lambda=\{\sum_{k=0}^{n-1}p^kl_k\,|\,l_k\in L,n\in\bn\}.$$
\end{itemize}
\end{enumerate}
\end{theorem}


\begin{proof}
\par
(i) First note that $m_B(x)=\frac12(1+e^{2\pi ix})$. Therefore the zeroes of $m_B$ are of the form $\frac k2$ with $k$ odd.
\par
Assume that $\mu_B$ is spectral and, by contradiction, that $r$ is odd. 
By a translation we can assume $0$ is in the spectrum $\Lambda$ of $\mu_{B}$. Take $\lambda,\lambda'\in\Lambda$, $\lambda\neq\lambda'$. We have then $\hat\mu_{B}(\lambda)=\hat\mu_{B}(\lambda')=\hat\mu_{B}(\lambda-\lambda')=0$, therefore, $m_B(p^{-n}(\lambda-\lambda'))=0$ for some $n\geq1$.
\par
Then for some $n_1,n_2,n_3\geq1$, $k_1,k_2,k_3$ odd, we have
$$\lambda=r^{n_1}\frac{k_1}{2a},\lambda'=r^{n_2}\frac{k_2}{2a},\lambda-\lambda'=r^{n_3}\frac{k_3}{2a}$$
But then if $n_1\geq n_2$, 
$$\lambda-\lambda'=\frac{r^{n_2}(r^{n_1-n_2}k_1-k_2)}{2a},$$
and $r^{n_1-n_2}k_1-k_2$ is even. This is a contradiction. So $r$ has to be even, and when $r$ is odd no $3$ exponential functions are mutually othogonal. 
\par
When $p$ is even $p=2a$, with $L=\{0,a\}$, it is easy to check that $(R,B,L)$ forms a Hadamard triple. Since $m_B$ has finitely many zeroes in any compact interval, we can apply the results in \cite{DuJo05}.
\par
If $p=2$ then $\mu_B$ is the Lebesgue measure on $[0,1]$. Since $W_B$ is $1$ only when $x\in\bz$, one can see easily that $\{0\}$ and $\{1\}$ are the $W_B$-cycles. Each of them contributes to the spectrum: the contribution of $0$ is 
$\{\sum_{k=0}^n2^ka_k\,|\,a_k\in\{0,1\},n\in\bn\}\bz_+\cup\{0\}$; the contribution of $1$ is $\{-2^n+\sum_{k=0}^{n-1}2^ka_k\,|\,a_k\in\{0,1\},n\in\bn\}=\bz_-$. Their union is the well-known spectrum $\bz$.
\par
If $p>2$, then $X_L$ is contained in $[0,\frac p2\sum_{k=1}^\infty p^{-k}]\subset[0,1)$ so the only $W_B$-cycle is $\{0\}$. With \cite{DuJo05}, we obtain the spectrum of $\mu_B$.
\par
(ii) We have
$$m_B(x_1,x_2)=\frac{1}{3}(1+e^{2\pi i x_1}+e^{2\pi ix_2}),\quad((x_1,x_2)\in\br^2).$$
If $(u,v)$ is a zero for $m_B$, then 
$$1=|e^{2\pi i v}|^2=|-1-e^{2\pi i u}|^2=2+2\cos(2\pi u),$$
so $\cos(2\pi u)=-1/2$, so $u\in \frac{1}{3}+\bz$ or $u\in\frac{2}{3}+\bz$. And, since $e^{2\pi iv}=-1-e^{2\pi iu}$,
$$(u,v)\in\left(\frac{1}{3},\frac{2}{3}\right)+\bz^2,\mbox{ or }(u,v)\in\left(\frac23,\frac13\right)+\bz^2.$$
\par
If $p$ is not a multiple of $3$ we see that the set of zeroes of $m_B$ in $\br^2$ is invariant under $R^T$, so we can use Theorem \ref{thnoort} to conclude that there no more than $3$ mutually orthogonal exponential functions in $L^2(\mu_B)$.
\par
If $p$ is a multiple of $3$, then take $L$ as stated in the theorem. Note that the attractor $X_L$ is contained in the segment $\{(x,-x)\,|\,x\in[-a,a]\}$ where $a=\frac{2p}3\sum_{k=1}^\infty p^{-k}<1$. With Theorem \ref{propwb}, we conclude that the only $W_B$-cycle is $\{0,0\}$. Also $m_B$ has finitely many zeroes in $X_L$ so we can use \cite{DuJo05} to conclude that the spectrum of $\mu_B$ is the one given in the statement of the theorem.
\par
(iii)
For this example, 
$$m_B(x,y,z)=\frac{1}{4}(1+e^{2\pi ix}+e^{2\pi iy}+e^{2\pi iz}).$$
We look for the zeroes of $W_B$. Let $(x,y,z)$ be a zero of $W_B$. Let $z_1=e^{2\pi ix}$, $z_1=e^{2\pi iy}$, $z_1=e^{2\pi iz}$. We have$$1+z_1+z_2+z_3=0.$$
If $z_1=-1$ then the other two are $z$ and $-z$. Suppose $z_1\neq -1$. Then $1+z_1=-(z_2+z_3)$. But $1+z_1$ is the diagonal of the rhombus formed  with $1$ and $z_1$. Similarly  for $z_2+z_3$. Since the two rhombbi have equal sides and equal diagonals, it follows that one of the sides $z_2$ or $z_3$ is $-1$, and the other is $-z_1$.
\par
Therefore one of the $z_i$'s is $-1$ and the other two form a pair $\{z,-z\}$.
\par
Thus we have that $(x,y,z)$ has one of the following forms
$(\frac{k_1}{2},a,a+\frac{k_2}2)$, $(a,\frac{k_1}{2},a+\frac{k_2}2)$, $(a,a+\frac{k_2}2,\frac{k_1}2)$, where $a\in\br$, and $k_1,k_2$ odd integers.
\par
Note that if $p$ is odd, the map $x\mapsto px\mod\bz^d$ leaves these forms invariant. With the notations in Theorem \ref{thnoort}, we have that the orbit $\mathcal{O}(Z)$ of the set of zeroes of $m_B$ inside $[0,1)^3$ has distance to $\bz^d$ at least $\frac{1}{2}$, because one of the components is always $\frac12$. We can use Theorem \ref{thnoort}, to obtain that we can have at most $(\left\lfloor\frac{\sqrt{3}}{1/2}\right\rfloor+1)^3=4^3=256$ mutually orthogonal exponential functions in $L^2(\mu_B)$.
\par
If $p$ is even, a simple computation shows that $(R,B,L)$ is a Hadamard triple. 
\par
When $p=2$, since $X_L$ is contained in $[0,1]^3$, using Theorem \ref{propwb} we get the the possible $W_B$-cycles are contained in the corners of this cube. By inspection we see that the points $(0,0,0), (1,1,0), (0,1,1), (1,0,1)$ are the $W_B$-cycles, and each of them forms a cycle of length $1$.
\par
Next, we look for the possible translations of invariant subspaces $V$ that contain minimal compact invariant sets (see \cite[Theorem 2.17]{DuJo06}. By \cite[Lemma 5.1]{DuJo06}, we must have $W_B((x,y,z)+v)=0$, for some point $(x,y,z)\in X_L$, and all $v\in V$, with $V$ invariant for $R^T$ and $(x,y,z)=\tau_l(x',y',z')$ for some $(x',y',z')\in X_L$ and some $l\in L$.
\par
Since $(x,y,z)$ is a zero for $W_B$ it has one of the particular forms described above. By symmetry, we can consider only the case when $(x,y,z)=(\frac{k_1}2,a,a+\frac{k_2}2)$.
\par
Since we want $(x,y,z)\in X_L$, this implies that $k_1=1$, $k_2=\pm 1$, and $a\in[0,1]$.
\par
This shows that the only subspaces $V$ that might appear here, are $\{(a,a,0)\,|\,a\in\br\}$, $\{(0,a,a)\,|\,a\in\br\}$ or $\{(a,0,a)\,|\,a\in\br\}$. 
Take $V$ to be the first one, the other subspaces can be treated identically, by symmetry. We need to have 
$W_B(\tau_l(x',y',z')+(a,a,0))=0$ for all $a\in\br$. Thus $\tau_l(x',y',z')$ has the form $(b,b\pm\frac12,\frac12)$. 
Since $(x',y',z')\in[0,1]^3$, it follows that $z'=0$, $x'=y'$ and $l=(1,0,1)$ or $l=(0,1,1)$.
\par
So, the translate of $V$ which is invariant is $0+V:=V=\{(a,a,0)\,|\,a\in\br\}$. Indeed, $\tau_l(V)\subset V$ for 
$l\in\{(0,0,0),(1,1,0)\}$, and $W_B(\tau_l(v))=0$ for all $v\in V$, if $l\in\{(1,0,1),(0,1,1)\}$. Thus $0+V$ is an invariant set.
\par
However, we claim that we can discard this invariant set. We restrict the maps $\tau_l$ to $0+V$. Since for $l\in\{(1,0,1),(0,1,1)\}$ we have $W_B(\tau_l(v))=0$ for all $v\in V$, we can discard those maps $\tau_l$. We have then an IFS with two maps $\tau_{(0,0,0)}$ and $\tau_{(1,1,0)}$ on $V$. We also have
$$W_B(\tau_{(0,0,0}v)+W_B(\tau_{(1,1,0)}v)=1,\quad(v\in V)$$
But $V$ is one-dimensional so $W_B$ has only finitely many zeroes in $V\cap X_L$ (actually, only one $(\frac12,\frac12,0)$). Therefore (see \cite{DuJo05}), the minimal compact invariant sets are the $W_B$-cycles inside $V$, which we already considered. Thus, there are no extra minimal compact invariant sets inside $0+V$, other than the $W_B$-cycles, and therefore the $W_B$-cycles form a complete list of minimal compact invariant sets.
\par
The contributions of each $W_B$-cycle to the spectrum is as described in the theorem (see \cite{DuJo05} and \cite[Lemma 4.9]{DuJo06}).
\par
If $p$ is even and $p>2$, then $X_L$ is contained in the cube $[0,a]^3$, with $a=\sum_{k=1}^\infty p^{-k}\frac p2=\frac{p}{2p-2}\leq\frac23$. Using Theorem \ref{propwb}, the only $W_B$-cycle is $\{0\}$. The invariant subspaces can be discarded as in the case $p=2$. And using \cite{DuJo06} we obtain the spectrum of $\mu_B$.
\end{proof}


\begin{lemma}\label{lemconf}
Let $p\in\bn$, $p\geq 2$, $R=pI_d$ and $B=\{0,e_1,\dots,e_d\}$. Denote by 
$$D_n:=\inf\left\{\left|1+\sum_{l=1}^de^{2\pi i\frac{k_l}{p^n}}\right|\,|\,k_i\in\bz\right\}.$$
If $\inf_{n\in\bn}p^nD_n>0$ then the measure $\mu_B$ is not spectral.
\end{lemma}

\begin{proof}
We check condition (ii) in Theorem \ref{thnoort}. Suppose this condition is not satisfied. Then for all $\epsilon>0$, there exists $n\in\bn$ and $x_1,\dots,x_d\in\br$ such that $m_B(x_1,\dots,x_d)=0$, and there exist $k_1,\dots,k_d\in\bz$ such that 
$$\left|p^nx_i-k_i\right|<\epsilon.$$
Then $|x_i-\frac{k_i}{p^n}|<\frac{\epsilon}{p^n}$ for all $i$. 
Since $m_B$ is a Lipschitz function, there is $M>0$ such that
$$|m_B(y_1,\dots,y_d)-m_B(y_1',\dots,y_d')|\leq M\max_{i}|y_i-y_i'|.$$
Then
$$|m_B(x_1,\dots,x_d)-m_B(\frac{k_1}{p^n},\dots,\frac{k_d}{p^n})|\leq M\frac{\epsilon}{p^n}.$$
Since $m_B(x_1,\dots,x_d)=0$ this implies 
$$\left|1+\sum_{l=1}^de^{2\pi i\frac{k_l}{p^n}}\right|\leq NM\frac{\epsilon}{p^n}.$$
Thus $D_n\leq NM\epsilon$. Letting $\epsilon$ tend to $0$, we obtain a contradiction.
\end{proof}


\begin{proposition}\label{propdiv}
Let $p\in\bn$, $p\geq 2$, $R=pI_d$ and $B=\{0,e_1,...,e_d\}$. Suppose $d+1$ is a linear combination, with non-negative integer coefficients, of some proper divisors of $p$. Then there exists and infinite family of mutually orthogonal exponential functions in $L^2(\mu_B)$. 
\end{proposition}

\begin{proof}
Suppose 
$$p=q_1d_1+...+q_sd_s$$
with $q_1,...,q_s\in\bn$ and $d_1,...,d_s$ proper divisors of $p$.
Define the point in $\br^d$, by repeating the $d_k$ roots of unity $q_k$ times, except for $d_1$ where we omit a $0$:
$$z_0=(\frac{1}{d_1},...,\frac{d_1-1}{d_1},\underbrace{\frac{0}{d_1},\frac{1}{d_1},...,\frac{d_1-1}{d_1}}_{q_1-1\mbox{ times}},...,\underbrace{\frac{0}{d_s},\frac{1}{d_s},...,\frac{d_s-1}{d_s}}_{q_s\mbox{ times}})$$
Then note that 
$$m_B(z_0)=\sum_{k=1}^sq_k\sum_{l=0}^{d_k-1}e^{2\pi i\frac{l}{d_k}}=0.$$
We prove that $$\{e_{p^nz_0}\,|\,n\geq 1\}$$ forms an orthogonal family.
Take $n>m\geq 1$. Then
$$p^{-m}(p^nz_0-p^mz_0)=p^{n-m}z_0-z_0\in\bz^d-z_0,$$
and this implies that $m_B(p^{-m}(p^nz_0-p^mz_0))=0$, and therefore $\hat\mu_B(p^nz_0-p^mz_0)=0$, so $e_{p^nz_0}$ and $e_{p^mz_0}$ are orthogonal. Since $n,m$ are arbitary this proves the proposition.
\end{proof}

\begin{corollary}
If $p$ is a multiple of $6$, then there is an infinite family of mutually orthogonal exponential functions in $L^2(\mu_B)$ no matter what $d$ is. 
\end{corollary}
\begin{proof}
The number $p$ has divisors $2$ and $3$. If $d+1$ is odd, then we can write $d+1=k\cdot 2+1\cdot 3$, and if $d+1$ is even then we can write $d+1=k\cdot 2+0\cdot 3$.
\end{proof}

\begin{remark}
In Lemma \ref{lemconf} and Proposition \ref{propdiv} we looked for zeroes of $m_B$ which have $p$-adic components. Applying $m_B$ to such zeroes, we obtain a vanishing sum of roots of unity. Such sums have been analyzed in a different cotext, in number theory. An interesting result which relates to our analysis is the following (see \cite{LaLe00}):
\par
If $n$ and $p$ are positive integers then there exists a vanishing sum of $n$ (not necessarily distinct) roots of unity of order $p$ if and only if $m$ can be written as a linear combination with non-negative integer coefficients of some proper divisors of $p$.
\par
Thus, we can find a zero $z_0$ of $m_B$ such that $p^nz_0\in\bz^d$ for some $n\geq1$ if and only if $d+1$ is a linear combination, with non-negative integer coefficients, of some proper divisors of $p$. So the hypotheses of Proposition \ref{propdiv} are optimal, if a similar argument is used.
\end{remark}
\begin{theorem}\label{thpmuld}
Let $p\in\bn$, $p\geq 2$, $R=pI_d$ and $B=\{0,e_1,...,e_d\}$. If $p$ is divisible by $d+1$ then $\mu_B$ is a spectral measure.
\begin{enumerate}
\item If $p=d+1$, then a spectrum for $\mu_B$ is $\bz (1,2,...,d)$.
\item If $p=m(d+1)$ with $m>1$, then a spectrum for $\mu_B$ is
$$\Lambda:=\{\sum_{k=0}^np^ka_k(m,2m,...,dm)\,|\,a_k\in\{0,...,d\},n\in\bn\}.$$

\end{enumerate}
\end{theorem}

\begin{proof}
We have $p=m(d+1)$ for some integer $m\geq1$. Let $v_0:=(1,2,...,d)\in\br^d$. Let 
$$L:=\{jmv_0\,|\,j\in\{0,1,...,d\}\}.$$
The matrix $(e^{2\pi i R^{-1}b\cdot l})_{b\in B,l\in L}$ is $(e^{2\pi i \frac{kj}{d+1}})_{k,j\in\{0,..,d\}}$, which is the matrix of the Fourier transform on the group $\bz_{d+1}$, so it is unitary. Thus $(R,B,L)$ forms a Hadamard pair. 
\par
Note also that the elements in $L$ are all multiples of $v_0$. Since $R$ leaves the subspace $v_0\br$ invariant, we can easily see that the segment $[0,\frac{d}{p-1}]$ is invariant for all maps $\tau_l$, $l\in L$. Therefore $X_L$ is contained in this segment. 
\par
The restriction of $m_B$ to any segment has finitely many zeroes, therefore $m_B$ has finitely many zeroes in $X_L$. We can then apply the results in \cite{DuJo05} to conclude that $\mu_B$ is spectral. 
\par
If $p=d+1$, then with Theorem \ref{propwb}, we see that there are two $W_B$-cycles, $\{0\}$ and $\{v_0\}$. The contributions of these $W_B$-cycles are
$$\Lambda(0):=\{\sum_{k=0}^n(d+1)^ka_kv_0\,|\,a_k\in\{0,...,d\},n\in\bn\}$$
$$\Lambda(v_0):=\{-(d+1)^nv_0+\sum_{k=0}^n(d+1)^ka_kv_0\,|\,a_k\in\{0,...,d\},n\in\bn\}.$$
Their union is $\bz v_0$ and it is a spectrum for $\mu_B$.
\par
If $p=m(d+1)$ with $m>1$, then with Theorem \ref{propwb}, $\{0\}$ is the only cycle so the spectrum of $\mu_B$ is as in (ii).
\end{proof}

\section{Concluding remarks}
We have studied a class of irregular patterns of points in $\br^d$
arising from finite families of affine mappings given by a fixed scaling
matrix $R$, and a fixed set of vectors $B$ in $\br^d$. 
      We show that iterations in the small leads to fractal sets and fractal
measures $\mu_B$. Motivated by analogies to lacunary Fourier series \cite{Kah86} we
further use iterations in the large to construct orthogonal complex
exponentials in the Hilbert spaces $L^2(\mu_B)$. 
      While in sections 3 and 4 we present general results, in section 5 we
restrict the discussion to a class of Sierpinski examples in $\br^d$. The
analysis of these examples shows that there is a certain rigidity which
limits the possibilities for $d$ small; while if $d = 4$ or higher, there some
unexpected additional possibilities.
     Our analysis of these higher dimensional cases suggests the following:
\begin{conjecture}
     Let $R$ be an expansive $d$ by $d$ matrix over $\bz$, and set $S = R^T$ the
transposed matrix. Let $B$ and $L$ be two subsets of $\br^d$ of the same cardinality
$N$, and assume that the Hadamard axiom in Definition \ref{defhada} holds. Let $\mu_B$ and
$\mu_L$ be the two Hutchinson measures booth with weights $p_i = 1/N$
corresponding to the respective sides in the dual pair of AIFSs $(R, B)$ and
$(S, L)$.
\begin{itemize}
\item[(a)]  Then  $\mu_B$ is a spectral measure if and only if $\mu_L$ is.
\item[(b)]  The two subsets $\Lambda_B$ and $\Lambda_L$  in $\br^d$  that index the complex
exponentials making up the ONBs for the Hilbert spaces $L^2(\mu_B)$ and
$L^2(\mu_L)$ are related by an explicit formula.
\end{itemize}
\end{conjecture}
\begin{acknowledgements}
The authors thank Brian Treadway for expert help with \TeX, 
and with programming and production of the graphics used in the figures. The authors wish to thank an anonymous referee for making helpful suggestions regarding the presentation of our results.

\end{acknowledgements}

\end{document}